\documentclass[aop,MSNbibl,seceqn,dvips]{arximspdf}
\usepackage{mathbh}
\usepackage{accents}
\usepackage{graphicx}

\doi{10.1214/12-AOP762} 
\volume{41}
\issue{3A}
\pubyear{2013}
\firstpage{1115}
\lastpage{1159}

\makeatletter

\renewcommand{\underline}{\underbar}
\renewcommand{\underbar}{\underaccent{\bar}}

\newcommand{\rrVert}{\Vert}
\newcommand{\rrvert}{\vert}
\newcommand{\llVert}{\Vert}
\newcommand{\llvert}{\vert}

\newtheorem{itlemma}{Lemma}[section]
\newtheorem{itproposition}[itlemma]{Proposition}
\newtheorem{theorem}[itlemma]{Theorem}
\newtheorem{itcorollary}[itlemma]{Corollary}

\newproclaim{itremark}[itlemma]{Remark}

\newcommand{\cal}{\mathcal}

\newcommand{\T}{{\cal{T}}}
\newcommand{\I}{{\mathsf{I}}}
\newcommand{\A}{{\cal{A}}}
\renewcommand{\H}{{\cal{H}}}
\newcommand{\GG}{{\cal{G}}}
\newcommand{\YY}{{\cal Y}}
\newcommand{\C}{{\cal{C}}}
\newcommand{\D}{{\cal{D}}}
\renewcommand{\S}{{\cal{S}}}

\newcommand{\B}{{\mathbb B}}
\newcommand{\N}{{\mathbb N}}
\newcommand{\Z}{{\mathbb Z}}
\newcommand{\R}{{\mathbb R}}
\renewcommand{\P}{{\mathbb P}}
\newcommand{\E}{{\mathbb E}}
\newcommand{\id}{{\mathbf{1}}}

\newcommand{\bs}{\setminus}

\newcommand{\NN}{\mathcal N}
\newcommand{\var}{\operatorname{var}}
\renewcommand{\AA}{{\mathbb A}}
\newcommand{\ind}{{\mathbh1}}

\makeatother

\begin{document}
\begin{frontmatter}

\title{From logarithmic to subdiffusive polynomial fluctuations
for internal DLA and related growth~models\thanksref{T1}}
\runtitle{Fluctuations for internal DLA}

\thankstext{T1}{Supported by GDRE 224 GREFI-MEFI,
the French Ministry of Education through the ANR BLAN07-2184264 grant,
by the European Research Council through the ``Advanced
Grant'' PTRELSS 228032.}

\begin{aug}
\author[A]{\fnms{Amine} \snm{Asselah}\ead[label=e1]{amine.asselah@univ-paris12.fr}}
\and
\author[B]{\fnms{Alexandre} \snm{Gaudilli\`ere}\corref{}\ead[label=e2]{gaudilli@cmi.univ-mrs.fr}}
\runauthor{A. Asselah and A. Gaudilli\`ere}
\affiliation{Universit\'e Paris-Est and Universit\'e de Provence}
\dedicated{Dedicated to Joel Lebowitz, for his 80th birthday}
\address[A]{LAMA\\
Universit\'e Paris-Est\\
61 avenue du g\'en\'eral de Gaulle\\
94010 Cr\'eteil Cedex\\
France\\
\printead{e1}} 
\address[B]{LATP\\
Universit\'e de Provence\\
CNRS, 39 rue F. Joliot Curie\\
13013 Marseille\\
France\\
\printead{e2}}
\end{aug}

\received{\smonth{5} \syear{2010}}
\revised{\smonth{4} \syear{2012}}

%
\begin{abstract}
We consider a cluster growth model on $\Z^d$, called internal diffusion
limited aggregation (internal DLA). In this model, random walks start
at the origin, one at a time, and stop moving when reaching a site not
occupied by previous walks. It is known that the asymptotic shape of
the cluster is spherical. When dimension is 2 or more, we prove that
fluctuations with respect to a sphere are at most a power of the
logarithm of its radius in dimension $d\ge2$. In so doing, we introduce
a closely related cluster growth model, that we call \textit{the
flashing process}, whose fluctuations are controlled easily and
accurately. This process is coupled to internal DLA to yield the
desired bound. Part of our proof adapts the approach of Lawler, Bramson
and Griffeath, on another space scale, and uses a sharp estimate
(written by Blach\`ere in our Appendix) on the expected
time spent by a
random walk inside an annulus.
\end{abstract}

%

\end{frontmatter}

\section{Introduction}\label{intro}

The internal DLA cluster of volume $N$, say $A(N)$,
is obtained inductively as follows. Initially, we assume
that the explored region is empty, that is, $A(0)=\varnothing$.
Then, consider $N$ independent discrete-time
random walks $S_1,\ldots,S_N$ starting from 0.
For $k\le N$, assume $A(k-1)$ is obtained, and define
\[
\tau_k=\inf\bigl\{ {t\ge0\dvtx S_k(t)\notin A(k-1)}
\bigr\} \quad\mbox{and}\quad A(k)=A(k-1)\cup\bigl\{S_k(\tau_k)
\bigr\}.
\]
In such a particle system, we call explorers the particles.
We say that the $k$th explorer is \textit{settled} on $S_k(\tau_k)$
after time
$\tau_k$, and is \textit{unsettled} before time~$\tau_k$. The cluster $A(N)$
consists of the positions of the $N$ settled explorers.

The mathematical model of
internal DLA was introduced first in the chemical physics literature
by Meakin and Deutch~\cite{meakin-deutch}.
There are many industrial processes that look like
internal DLA; see the nice review paper~\cite{landolt}.
The most important seems to be electropolishing,
defined as \textit{the improvement of surface finish of a metal effected
by making it anodic in an appropriate solution}. There are actually
two distinct industrial processes (i) \textit{anodic leveling or smoothing}
which corresponds to the elimination of surface roughness of
height larger than 1 micron, and (ii) \textit{anodic brightening} which
refers to elimination of surface defects which are protruding by
less than 1 micron. The latter phenomenon
requires an understanding of atom removal from a crystal lattice.
It was noted in~\cite{meakin-deutch} that, at a qualitative level,
the model produces smooth clusters, and the authors wrote,
``it is also of some fundamental significance to know just how smooth
a surface formed by diffusion limited processes may be.''

Diaconis and Fulton~\cite{diaconis-fulton} introduced internal DLA
in mathematics. They allowed explorers to start on distinct sites, and
showed that the law of the cluster was invariant under permutation of
the order in which explorers were launched. This invariance,
named \textit{the abelian property}, was central in their motivation.
They treat, among other things, the special one-dimensional case.

In dimension two or more,
Lawler, Bramson and Griffeath~\cite{lawler92} prove that in order
to cover, without holes, a sphere of radius $n$, we need about the
number of sites of $\Z^d$ contained in this sphere. In other words,
the asymptotic shape of the cluster
is a sphere. Then, Lawler in~\cite{lawler95} shows
subdiffusive fluctuations.
The latter result is formulated in terms of inner and outer errors,
which we now introduce with some notation.
We denote with $\|\cdot\|$ the Euclidean norm on $\R^d$.
For any $x$ in $\R^d$ and $r$ in $\R$, set
\[
B(x,r) = \bigl\{ y\in\R^d\dvtx\|y-x\| < r \bigr\}
\quad\mbox{and}\quad\B(x,r) =
B(x,r) \cap\Z^d.
\]
For $\Lambda\subset\Z^d$,
$|\Lambda|$ denotes the number of sites in $\Lambda$.
The inner error $\delta_I(n)$ is such that
\[
n-\delta_I(n)=\sup\bigl\{ {r \geq0\dvtx\B(0,r)\subset A\bigl(\bigl|\B
(0,n)\bigr|\bigr)} \bigr\}.
\]
Also, the outer error $\delta_O(n)$ is such that
\[
n+\delta_O(n)=\inf\bigl\{ {r \geq0\dvtx A\bigl(\bigl|\B(0,n)\bigr|\bigr)
\subset\B(0,r)} \bigr\}.
\]
The main result of~\cite{lawler95} reads as follows.
%
%
\begin{theorem}[(Lawler)]\label{theo-lawler}
Assume $d\ge2$. Then
%
%
\begin{equation}
\label{lawler-results} P \bigl( {\exists n(\omega)\dvtx\forall n\ge
n(\omega)\
\delta_I(n)\le n^{1/3}\log(n)^2} \bigr)=1\vadjust{\goodbreak}
\end{equation}
and
%
%
\begin{equation}
\label{lawler-outer} P \bigl( {\exists n(\omega)\dvtx\forall n\ge
n(\omega)\
\delta_O(n)\le n^{1/3}\log(n)^4} \bigr)=1.
\end{equation}
\end{theorem}
%
Since Lawler's paper, published 15 years ago,
no improvement of these estimates was achi\-eved,
but it is believed that fluctuations are on a much smaller scale
than $n^{1/3}$. Moreover, (\ref{lawler-results}) and
(\ref{lawler-outer}) are almost sure
upper bounds on errors, and no lower bound on the inner or outer error
has been established. Computer simulations~\cite{machta-moore,levine}
suggest indeed that fluctuations are logarithmic.
In addition, Levine and Peres studied
a deterministic analogue of internal DLA, the rotor-router model,
introduced by Propp~\cite{kleber-propp}.
They bound, in~\cite{levine-peres}, the inner error $\delta_I(n)$
by $\log(n)$, and the outer error $\delta_O(n)$ by $n^{1-1/d}$.

Our main result is the following improvement of Theorem~\ref{theo-lawler}.
%
%
\begin{theorem}\label{prop-main}
Assume $d\ge2$.
There is a positive constant $A_d$ such that
%
%
\begin{equation}
\label{ag-inner} P \bigl( {\exists n(\omega)\dvtx\forall n\ge n(\omega
)\
\delta_{I}(n)\le A_d\log(n) } \bigr)=1
\end{equation}
and
%
%
\begin{equation}
\label{ag-outer} P \bigl( {\exists n(\omega)\dvtx\forall n\ge n(\omega
)\
\delta_{O}(n)\le A_d \log^2(n) } \bigr)=1.
\end{equation}
\end{theorem}

\textit{Note added in proof.}
At about the same time,
and with an independent approach, Jerison, Levine and
Sheffield~\cite{JLS1} obtained similar results with
an improved bound on the outer error in $d=2$. Then, by refining
our approach, we obtained in~\cite{AG2}
a bound of order $\sqrt{ \log(n)}$ for
both internal and external errors in dimension three or more.
Jerison, Levine and Sheffield~\cite{JLS2}
did the same by following their approach.

Our approach builds on the work of Lawler, Bramson and Griffeath
\cite{lawler92}, which we review later. It also deals with
more general models of diffusion limited aggregation which we now
describe.
Indeed, we introduce a family of cluster growth models for which a
control of the fluctuations of the cluster shape is easily obtained.
These growth models are built so that the asymptotic shape is spherical,
but still they exhibit a large diversity of fluctuations parametrized
by a certain width ranging from a large constant to a power $1/3$ of the
radius of the asymptotic sphere.
Moreover, all these clusters are coupled to internal DLA,
and, as a consequence, we obtain logarithmic bounds
on the fluctuations for internal DLA.
We generalize internal DLA by allowing explorers to
settle only at some special times. Thus, each explorer $i$ is associated
with a collection of times $\{\sigma_{i,k}, k\in\N\}$ and
\[
\tau_i^*=\inf\bigl\{ {\sigma_{i,k}\dvtx S_i(
\sigma_{i,k})\notin A^*(i-1)} \bigr\}\quad\mbox{and}\quad
A^*(i)=A^*(i-1)\cup
\bigl\{S_i\bigl(\tau_i^*\bigr)\bigr\}.
\]
The internal DLA is recovered as we choose $\sigma_{i,k}=k$ for
all $i=1,\ldots,N$ and $k\in\N$. We call $\{\sigma_{i,k}, k\in\N\}$
the \textit{flashing times} associated to the $i$th explorer, and
$\{S_i(\sigma_{i,k}), k\in\N\}$ its \textit{flashing positions}.

\begin{figure}

\includegraphics{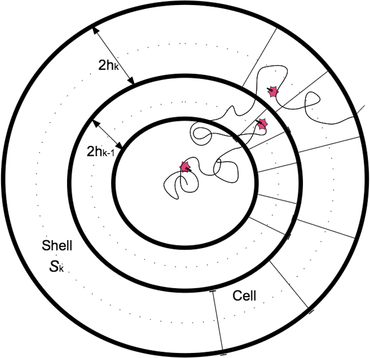}

\caption{Cell decomposition, and flashing positions as stars.}
\label{figflashing}
\end{figure}

In this paper, we consider stopping times of a special form,
linked with the spherical nature of the internal DLA cluster.
An illustration with one \textit{flashing} explorer's trajectory
is made in Figure~\ref{figflashing}.

%
%

The precise definition of the \textit{flashing times} requires
additional notation, which we postpone
to Section~\ref{sec-flashing}. We describe here
key features of flashing processes.
We first choose a sequence of widths, say $\H=\{h_n, n\in\N\}$, and
then partition $\Z^d$ into concentric
shells $\{\S_n, n\in\N\}$, whose respective widths
are $\{2h_n, n\in\N\}$.
Each shell is in turn partitioned into cells, which are brick-like domain,
of side length equal to the width of the shell. The flashing times are
chosen such that
(i) an explorer flashes at most once in each shell,
(ii) the flashing position, in a shell, is essentially uniform
over the cell an explorer first hits upon entering the shell and
(iii) when an explorer leaves a shell, it cannot afterward
flash in it.

For a given sequence $\H$, we call the process just described the
$\H$-flashing process.
Note that feature (ii) is the seed of a deep difference with internal DLA.
\textit{The mechanism of covering a cell, for the flashing process},
\textit{is very much the same
as completing an album in the classical coupon-collector process}.
Thus, we need of the order of $V\log(V)$ explorers
to cover a cell of volume $V$. For internal DLA, with explorers
started at the origin, we only need of order $V$ explorers to
cover a sphere of volume $V$ as shown in~\cite{lawler92}, and
we believe that we need a number of explorers
of order $|{\cal C}|$ to cover a cell ${\cal C}$, even if they start
on the boundary of the cell.
In addition, feature (ii)
allows us to \textit{localize} the covering mechanism, in the
sense that a particle entering a shell
cannot flash outside the cell through which it entered that shell.
Finally, feature (iii) is essential for having a useful coupling
between flashing and internal DLA processes.
%
%
\begin{itlemma}\label{theo-coupling}
Assume that $N$ is an integer, and
$\H$ is a sequence of positive integers.
There is a coupling between the two processes,
using the same trajectories $S_1,\ldots,S_N$ such that
%
%
\begin{equation}
\label{feature-coupling} A(N)=\bigcup_{i=1}^N
\bigl\{S_i\bigl(T(i)\bigr)\bigr\}\quad\mbox{and}\quad A^*(N)=\bigcup
_{i=1}^N \bigl\{S_i
\bigl(T^*(i)\bigr)\bigr\}
\end{equation}
and $T^*(i)\ge T(i)$ for all $i=1,\ldots,N$.
\end{itlemma}
As a corollary of Lemma~\ref{theo-coupling}, we have
the following useful result.
%
%
\begin{itcorollary}\label{cor-coupling}
Under the hypotheses of the previous lemma, for $k\geq1$:
\begin{itemize}
\item
if $A^*(N)\subset\bigcup_{j<k} \S_j$,
then $A(N)\subset\bigcup_{j<k} \S_j$;
\item
if $\bigcup_{j<k} \S_j\subset A^*(N)$,
then $\bigcup_{j<k} \S_j\subset A(N)$.
\end{itemize}
\end{itcorollary}

An $\H$-flashing process, with $h_j\ge h_0$ for $j\ge0$, and
$h_0$ a large constant, produces a cluster $A^*(N)$, for
which we bound easily the inner error, $\delta^*_I(n)$.
Then, to bound the outer error, $\delta^*_O(n)$,
we follow the approach of~\cite{lawler95},
though with a slightly simpler proof.
%
%
\begin{itproposition}\label{theo-flashing}
Assume that for $j\ge1$, $h_j\le h_{j+1}\le(1+\frac{1}{2j})h_j$,
with a large $h_0$.
For a positive constant $A_d^*$, we have
%
%
\begin{equation}
\label{results-flashing} P \bigl( {\exists n(\omega)\dvtx\forall n\ge
n(\omega)\
\delta^*_{I}(n)\le A_d^*h(n)\log(n) } \bigr)=1
\end{equation}
and
%
%
\begin{equation}
\label{results-outer} P \bigl( {\exists n(\omega)\dvtx\forall n\ge
n(\omega)\
\delta^*_{O}(n)\le A_d^*h(n)\log^2(n) }
\bigr)=1,
\end{equation}
where $h(n) = \max\{h_k\in{\mathbb R}\dvtx r_k\leq n\}$.
\end{itproposition}

Finally, we establish lower bound on the inner and outer error.
%
%
\begin{itproposition}\label{theo-optimal}
Assume that $h_0$ is large enough. Then, there is a constant $a_d^*$
such that
%
%
\begin{equation}
\label{lower-inner} P \bigl( {\exists n(\omega)\dvtx\forall n\ge
n(\omega)\
\delta^*_{I}(n)\ge a_d^* h(n)\log\bigl(h(n)\bigr)}
\bigr)=1
\end{equation}
and
%
%
\begin{equation}
\label{lower-outer-th} P \bigl( {\exists n(\omega)\dvtx\forall n\ge
n(\omega)\
\delta^*_{O}(n)\ge a_d^* h(n)\log(n)} \bigr)=1.
\end{equation}
\end{itproposition}
Corollary~\ref{cor-coupling} and Proposition~\ref{theo-flashing}, with
the choice $h_j=h_0$ for all $j>0$, imply Theorem~\ref{prop-main} which
deals with internal DLA.

Let us now review previous work on internal DLA.

\subsection*{On previous bounds for internal DLA}
We describe the approach of
\cite{lawler92}, for establishing the upper bound for the
inner error. It is convenient to consider explorers starting
outside the origin with initial configuration denoted~$\eta$.
We denote also by $A(\Lambda,\eta)$ the cluster obtained from explorers
initially on~$\eta$, with an explored region $\Lambda\subset\Z^d$.

Now, for a site $z\in\Z^d$, we call $W(\eta,z)$ [resp., $M(\eta,z)$]
the number of explorers (resp., of random walks) which visit $z$
before settling. For an integer~$n$, and $\eta$ consisting of
$|\B(0,n)|$ explorers at the origin, the authors
of~\cite{lawler92} first write
\[
\bigl\{ {\B(0,r)\not\subset A(\varnothing, \eta)} \bigr\} \subset\bigcup
_{z\in\B(0,r)} \bigl\{ {W(\eta,z)=0} \bigr\}.
\]
Then, they look for the largest value of $r_n$ (in terms of $n$) which
guarantees that
$| \B(0,r_n)|\times\sup_{z\in\B(0,r_n)}P(W(\eta,z)=0)$ be the
term\vspace*{1pt} of
a convergent series.

The approach of~\cite{lawler92} is based on the following observations.
(i) If explorers would not settle,
they would just be independent random walks; (ii) exactly one
explorer occupies each site of the cluster. Thus, the following
equality holds in law:
\[
W(\eta,z)+M\bigl(A(\varnothing,\eta),z\bigr) \ge M(\eta,z).
\]
Now, an observation of Diaconis and Fulton~\cite{diaconis-fulton} is that
we can realize the cluster by sending many \textit{exploration waves}.
Let us illustrate this observation with two waves. We first stop the
explorers on the external
boundary of a ball of radius~$R$, say $\partial\B(0,R)$.
The cluster consisting of the positions of settled explorers is denoted
$A_R(\varnothing,\eta)$, so that
$A_R(\varnothing,\eta)\subset\B(0,R)$.
The configuration with stopped explorers on $\partial\B(0,R)$
is denoted $\zeta_R(\eta)$.
%
Then,
the second wave consists in launching the explorers of
$\zeta_R(\eta)$, with explored region
$A_R(\varnothing,\eta)$. In other words,
we have an equality in law
\[
A(\varnothing,\eta)=A_R(\varnothing,\eta)\cup A \bigl(
{A_R(\varnothing,\eta),\zeta_R(\eta)} \bigr).
\]
Moreover, if the index $R$ refers only to explorers
(or walks) of the first wave, then for
$z\in\B(0,R)$,
%
%
\begin{equation}
\label{lawler-3} W_R(\eta,z)+M_R\bigl(A_R(
\varnothing,\eta),z\bigr) \ge M_R(\eta,z).
\end{equation}
The authors of~\cite{lawler92} consider
$R=n$ and $z\in\B(0,n)$. Since
$W(\eta,z)\ge W_n(\eta,z)$, we have using (\ref{lawler-3}),
for any $\alpha>0$,
%
%
\begin{equation}
\label{main-lawler92} P \bigl( {W ( {\eta,z} )=0} \bigr)\le P \bigl(
{M_n ( {\eta,z} )<\alpha} \bigr)+P \bigl( {M_n \bigl( {
\B(0,n),z} \bigr)>\alpha} \bigr).
\end{equation}
We then look for sites $z$ such that
$E[M_n(\eta,z)]>\alpha>E[M_n(\B(0,n),z)]$ (and $\eta=|\mathbb B(0,n)|\delta_0$).
Note that $M_n(\eta,z)$ and $M_n(\B(0,n),z)$ are sums of independent
Bernoulli variables with well-known large deviation estimates.
If we set
\[
2\alpha=E\bigl[M_n(\eta,z)\bigr]+E\bigl[M_n\bigl(
\B(0,n),z\bigr)\bigr]
\]
and
\[
\tilde\mu_n(z)=E
\bigl[M_n(\eta,z)\bigr]-E\bigl[M_n\bigl(\B(0,n),z\bigr)
\bigr],
\]
then
%
%
\begin{eqnarray}
\label{estimate-lawler} P \bigl( {M_n(\eta,z)<\alpha} \bigr)&\le&\exp
\biggl( {-\frac
{(E[M_n(\eta,z)]-\alpha)^2}{ 2E[M_n(\eta,z)]}} \biggr)\nonumber\\
[-8pt]\\[-8pt]
&\le&\exp\biggl(
{-\frac{ \tilde\mu_n^2(z)}{8 E[M_n(\eta,z)]}} \biggr).\nonumber
\end{eqnarray}
Lawler in~\cite{lawler95} establishes that for $z\in\B(0,n)$,
\[
E \bigl[ {M_n(\eta,z)} \bigr]\sim n\bigl(n-\|z\|\bigr)\quad
\mbox{and}\quad\tilde
\mu_n(z)\sim\bigl( {n-\|z\|} \bigr)^2.
\]
Replacing these values in (\ref{estimate-lawler}), the bound
$n-\|z\|\ge n^{1/3}\log(n)$ is such that
$P(W_n(\eta,z)=0)$ is the term of a convergent series.

We now sketch our
main ideas leading to logarithmic fluctuations for internal DLA.

\subsection*{On logarithmic fluctuations}
Our approach is inspired by Lawler, Bramson and Griffeath's work
\cite{lawler92}. We develop three original ideas:
(i) we propose a cluster growth model,
the \textit{flashing process}, whose covering mechanism
is simpler than internal DLA;
(ii) we look at an intermediary scale, \textit{the scale of cells},
since the deviations of the number of visits decrease with
the cell-length;
(iii) we build a coupling between \textit{flashing process} and
internal DLA
which allows us to transport bounds from one model to the other.

Let us describe how the idea of an intermediary scale is
used in the context of flashing processes. Recall that we first partition
$\Z^d$ into a sequence of concentric shells.
Each shell is partitioned into cells
whose side length equals the width of the shell.
Now, we observe that a site has good chances
to lie inside the cluster if some cell, say $\C$, about this site,
is crossed by \textit{many} explorers.
The notation $W(\eta,\C)$ refers to the number
of explorers visiting~$\C$, when their initial configuration is $\eta$.
We drop the index $n$ appearing in $W_n(\eta,z)$ since there are no
more constraints on not escaping the ball $\B(0,n)$.
Now, the coupon-collector nature of the covering mechanism
suggests that for some
positive constant $\alpha_d$,
%
%
\begin{eqnarray}
\label{size-C}
&&
W(\eta,\C)\ge\alpha_d
|\C|\times\log\bigl(|\C|\bigr)\nonumber\\[-8pt]\\[-8pt]
&&\qquad\Longrightarrow\quad\C\subset A(\varnothing,\eta) \mbox{ with a large
probability.}\nonumber
\end{eqnarray}
We neglect in these heuristics the $\log(|\C|)$ term in
(\ref{size-C}).

Note that in~\cite{lawler92}, all the explorers start from
the origin, whereas here, we only know that they \textit{cross $\C$}.
For internal DLA, estimating the probability that $\C$ is not
covered, when $\C$ is large and
$W(\eta,\C)\ge\alpha_d |\C|$ raises a difficulty which is absent
when considering flashing processes.

We now make our argument more precise.
For a scale $h$ and an integer $K>1$, to be determined,\vadjust{\goodbreak}
assume that $\B(0,n-Kh)$ is covered by settled explorers.
Partition the shell $\S=\B(0,n-(K-1)h)\bs\B(0,n-Kh)$
into about $(n/h)^{d-1}$ \textit{cells}, each of volume $h^d$.
It is also convenient to stop the explorers as they reach the boundary
of $\B(0,n-Kh)$. Thus, with such a stopped process, explorers are
either settled inside $\B(0,n-Kh)$ or
unsettled but stopped on its boundary,
denoted $\partial\B(0,n-Kh)$. What we have called earlier
\textit{the number of explorers crossing $\C$} is taken here
to be the unsettled explorers stopped on $\C\cap\partial\B(0,n-Kh)$.

Assuming (\ref{size-C}) holds, it remains to show that
the probability of the
event $\{\exists\C\in\S\dvtx W(\eta,\C)<\alpha_d |\C|\}$ is small.
We improve (\ref{main-lawler92}) by first using the independence
between $W(\eta,\C)$ and $M(\B(0,n-Kh),\C)$, and then by
replacing $A_R(\varnothing,\eta)$ by $\B(0,n-Kh)$ in (\ref{lawler-3})
with $R=n-Kh$ and $\eta=|\B(0,n)|\ind_0$,
%
%
\begin{equation}
\label{main-AG} W(\eta,\C)+M\bigl(\B(0,n-Kh),\C\bigr)\ge M(\eta,\C).
\end{equation}
Also, we define
\[
\mu(\C)=E\bigl[M(\eta,\C)\bigr]-E\bigl[M\bigl(\B(0,n-Kh),\C\bigr)\bigr].
\]
Now, using that $M(\eta,z)$ and
$M(\B(0,n-Kh),z))$ are sums of independent Bernoulli variables,
we show that (\ref{main-AG}) implies a Gaussian-type lower tail
%
%
\begin{equation}
\label{intro-amin1} P \bigl( {W(\eta,\C)<\alpha_d |\C|} \bigr)\le
\exp\biggl( {-\frac{ ( {\mu(\C)-\alpha_d |\C|}
)^2}{c\nu(\C)}} \biggr)
\end{equation}
for a positive constant $c$, and where $\nu(\C)$
\[
\nu(\C)=\var\bigl( {M(\eta,\C)} \bigr)-\var\bigl( {M\bigl(\B(0,n-K
h),\C
\bigr)} \bigr).
\]
We then show that both $\mu(\C)$ and $\nu(\C)$ are of order $K|\C|$.
Then, $P(W(\eta,{\cal C}) < \alpha_d|{\cal C}|)$
is summable as soon as $K|\C|\ge A\log(n)$.

\subsection*{Outline of the paper}
The rest of the paper is organized as follows.
Section~\ref{sec-notation} introduces the main notation, and
recalls known useful facts.
In Section~\ref{sec-flashing}, we build the flashing process,
give an alternative construction through exploration waves
and sketch the proof of Lemma~\ref{theo-coupling}.
In Section~\ref{sec-results}, we prove Propositions~\ref{theo-flashing}
and~\ref{theo-optimal} using the construction in terms of
exploration waves. In Section~\ref{sec-harmonic},
we obtain a sharp estimate on the expected number of explorers
crossing a given cell, and prove feature (ii) of the flashing times.
Both proofs are based on classical potential theory estimates.
Finally, in the \hyperref[app]{Appendix}, we give a proof of Lemma~\ref
{theo-coupling},
and recall a result of S\'ebastien Blach\`ere.

\section{Notation and useful tools}\label{sec-notation}
\subsection{Notation}
We say that $z,z'\in\Z^d$ are nearest neighbors when $\|z-z'\|=1$, and
we write $z\sim z'$. For any subset $\Lambda\subset\Z^d$, we define
\[
\partial\Lambda= \bigl\{ {z\in\Z^d\bs\Lambda\dvtx\exists
z'\in\Lambda, z'\sim z} \bigr\}.\vadjust{\goodbreak}
\]
For any $r\leq R$, we define the annulus
%
%
\begin{equation}
\label{annulus-dfn} A(r,R) = B(0,R)\setminus B(0,r)
\quad\mbox{and}\quad{\mathbb
A}(r,R) =
A(r,R)\cap\Z^d.
\end{equation}
A trajectory $S$ is a discrete nearest-neighbor path on $\Z^d$.
That is, $S\dvtx\N\to\Z^d$ with $S(t)\sim S(t+1)$ for all integer $t$.
For a subset $\Lambda$ in $\Z^d$, and a trajectory $S$,
we define the hitting time of $\Lambda$ as
\[
H(\Lambda;S)=\min\bigl\{t\geq0 \dvtx S(t)\in\Lambda\bigr\}.
\]
We often omit $S$ in the notation when no confusion is possible.
We use the shorthand notation
\[
B_n=B(0,n),\qquad\B_n=\B(0,n),\qquad H_R=H
\bigl(B_R^c\bigr)\quad\mbox{and}\quad H_z = H\bigl(
\{z\}\bigr).
\]
For any $a$, $b$ in $\R$ we write $a\wedge b = \min\{a,b\}$, and
$a\vee b=\max\{a,b\}$.
Let $\Gamma$ be a finite collection of trajectories
on $\Z^d$. For $R>0$, $z$ in $\Z^d$ and $\Lambda$ a subset
of~$\Z^d$, we call
$M(\Gamma, R,z)$ [resp., $M(\Gamma,R, \Lambda)$]
the number of trajectories which exit $\B(0,R)$ on $z$ (resp., in
$\Lambda$).
\[
M(\Gamma,R,z) = \sum_{S\in\Gamma} \id_{\{S(H_R)=z\}}
\quad\mbox{and}\quad M(\Gamma,R,\Lambda) = \sum_{z\in\Lambda} M(
\Gamma,R,z).
\]
When we deal with a collection of independent random trajectories, we
rather specify its initial configuration $\eta\in\N^{\Z^d}$,
so that $M(\eta,R,z)$ is the number of random walks
starting from $\eta$ and hitting $\B(0,R)^c$ on $z$. Two types of
initial configurations are important here:
(i) the configuration $n\id_{z^*}$ formed by $n$ walkers starting on
a given site $z^*$ and
(ii) for $\Lambda\subset\Z^d$,
the configuration $\id_\Lambda$
that we simply identify with $\Lambda$.
For any configuration $\eta\in\N^{\Z^d}$ we write
\[
|\eta| = \sum_{z\in\Z^d} \eta(z).
\]
For any $\Lambda\subset\Z^d$, we define Green's function
restricted to $\Lambda$, $G_\Lambda$, as follows. For $x,y\in\Lambda$,
the expectation with respect to
the law of the simple random walk started at $x$, is denoted with $\E_x$
(the law is denoted $\P_x$) and
\[
G_\Lambda(x,y)= \E_x \biggl[ \sum
_{0\le n < H(\Lambda^c)} \id_{\{S(n) = y\}} \biggr].
\]
In dimension 3 or more, Green's function on the whole space
is well defined and denoted $G$. That is,
for any $x,y\in\Z^d$,
\[
G(x,y)=\E_x \biggl[ \sum_{n\ge0}
\id_{\{S(n) = y\}} \biggr].
\]
In dimension 2, the potential kernel plays the role of
Green's function
\[
a(x,y)=\lim_{n\to\infty} \E_x \Biggl[ {\sum
_{l=0}^{n} \bigl(\id\bigl\{ {S(l)=x} \bigr\} -\id
\bigl\{ {S(l)=y } \bigr\} \bigr) } \Biggr].
\]
%
\subsection{Some useful tools}\label{sec-green}
We recall here some well-known facts.
Some of them are proved for the reader's convenience. This section
can be skipped at a first reading.

In~\cite{lawler92}, the authors emphasized the fact
that the spherical limiting shape of internal DLA
was intimately linked to strong isotropy properties
of Green's function. This isotropy is expressed
by the following asymptotics (Theorem~4.3.1 of~\cite{lawler-limic}).
In $d\ge3$, there is a constant $K_g$, such that for any $z\not=0$,
%
%
\begin{equation}
\label{main-green} \biggl|G(0,z) - \frac{C_d}{\|z\|^{d-2}} \biggr|
\le\frac{K_g}{\|
z\|^d} \qquad\mbox{with }
C_d=\frac{2}{v_d(d-2)},
\end{equation}
where $v_d$ stands for the volume of
the Euclidean unit ball in $\R^d$.
The first order expansion (\ref{main-green})
is proved in~\cite{lawler-limic}
for general symmetric walks
with finite $d+3$ moments and vanishing third
moment. All the estimates we use are eventually based on (\ref{main-green}),
and we emphasize the fact that the estimate is uniform in $\|z\|$.
There is a similar expansion for the potential kernel.
Theorem 4.4.4 of~\cite{lawler-limic} establishes that
for $z\not=0$ (with $\gamma$ the Euler constant),
%
%
\begin{equation}
\label{main-potential} \biggl|a(0,z) - \frac{2}{\pi}\log\bigl( {\|z\|}
\bigr)-
\frac{2\gamma+\log(8)}{\pi} \biggr|\le\frac{K_g}{\|z\|^2}.
\end{equation}
We recall a rough but useful result about the exit site distribution
from a sphere. This is Lemma 1.7.4 of~\cite{lawler}.
%
%
\begin{itlemma}\label{lem-lawler}
There are two positive constants $c_1,c_2$ such that for
any $z\in\partial B(0,n)$, and $n>0$
%
%
\begin{equation}
\label{hitting-position} \frac{c_1}{n^{d-1}}\le\P_0
\bigl(S(H_n)=z\bigr)\le\frac{c_2}{n^{d-1}}.
\end{equation}
\end{itlemma}
We now state an elementary lemma.
%
%
\begin{itlemma}\label{lem-border}
Each $z^*$ in $\Z^d\setminus\{0\}$
has a nearest-neighbor $z$ (i.e., $z^*\sim z$) such that
%
%
\begin{equation}
\label{goodnn} \|z\|\leq\bigl\|z^*\bigr\| -\frac{1}{2\sqrt{d}}.
\end{equation}
\end{itlemma}
\begin{pf}
Without loss of generality
we can assume that all the coordinates
of $z^*$ are nonnegative.
Let us denote by $b$ the maximum of these coordinates, and note that
%
%
\begin{equation}
\label{dfna} \bigl\|z^*\bigr\|^2 \leq db^2\quad\mbox{and}\quad b\ge1.
\end{equation}
Denote by $z$ the nearest-neighbor obtained from $z^*$
by decreasing by one unit a maximum coordinate.
Using (\ref{dfna}),
%
%
\begin{equation}
\label{simple} \bigl\|z^*\bigr\|^2 - \|z\|^2=b^2-(b-1)^2=2b-1
\ge b\ge\frac{\|z^*\|}{\sqrt{d}}.
\end{equation}
Note that (\ref{goodnn}) follows from $2\|z^*\|(\|z^*\| - \|z\|)\ge
\|z^*\|^2 - \|z\|^2$, and (\ref{simple}).\vadjust{\goodbreak}~%
\end{pf}

We state now a handy estimate dealing with sums of independent
Bernoulli variables.
%
%
\begin{itlemma}\label{lem-folk}
Let $\{X_n,Y_n, n\in\N\}$ be independent 0--1 Bernoulli variables.
For integers $n,m$ let $S=X_1+\cdots+X_n$ and
$S'=Y_1+\cdots+Y_m$. Define for $t\in\R$
\[
f(t)=e^t-1-t\quad\mbox{and}\quad g(t)= \bigl( {e^t-1}
\bigr)^2.
\]
If $0\le t\le\log(2)$, then
%
%
\begin{equation}
\label{pos-folk} \frac{E [ {\exp( {t ( {S-E[S]} )} )}
]} {
E [ {\exp( {t ( {S'-E[S']} )} )} ]} \le\exp\Biggl( {f(t)E \bigl[
{S-S'} \bigr]+g(t)\sum_{i=1}^m
E[Y_i]^2} \Biggr).
\end{equation}
Assume now that for $\kappa>1$, $\sup_n E[Y_n]\le\frac{\kappa-1} {
\kappa}$. If $t\le0$, then
%
%
\begin{equation}
\label{neg-folk}\quad\frac{E [ {\exp( {t ( {S-E[S]} )} )}
]} {
E [ {\exp( {t ( {S'-E[S']} )} )} ]} \le\exp\Biggl( {f(t)E \bigl[
{S-S'} \bigr]+\frac{\kappa}{2} g(t)\sum
_{i=1}^m E[Y_i]^2} \Biggr).
\end{equation}
\end{itlemma}
\begin{pf}
Let $X$ be a Bernoulli variable, and $p=E[X]$. Using the inequality
$e^x\ge1+x$ for $x\in\R$, we have
%
%
\begin{eqnarray}
\label{folk-1} E \bigl[ {\exp\bigl( {t\bigl(X-E[X]\bigr)} \bigr)}
\bigr]&=&pe^{t(1-p)}+(1-p)e^{-tp}\nonumber\\
&=&e^{-pt} \bigl( {1+p
\bigl(e^t-1\bigr)} \bigr) \\
&\le&\exp\bigl( {f(t) E[X]} \bigr).\nonumber
\end{eqnarray}
For a lower bound, we distinguish two cases.

First, assume $t\ge0$.
We claim that
$\exp(x-x^2)\le1+x$
for $0\le x\le1$.
Indeed, we use three obvious inequalities: $e^x\ge1+x$ for $x\in\R$,
(i) for $x\le1$, $1+x+x^2\ge e^x$, and (ii) $(1+x^2)(1+x)\ge1+x+x^2$.
Thus
\[
e^{x^2}(1+x)\ge\bigl(1+x^2\bigr) (1+x)\ge1+x+x^2
\ge e^x.
\]
This yields the claim. Now, set $x=p(e^t-1)$, so that $x\le1$ when
$e^t\le2$. The last inequality in (\ref{folk-1}) yields
%
%
\begin{eqnarray}
\label{folk-2} E \bigl[ {\exp\bigl( {t\bigl(X-E[X]\bigr)} \bigr)} \bigr]
&\ge&\exp
\bigl( {-tp+p\bigl(e^t-1\bigr)-p^2\bigl(e^t-1
\bigr)^2} \bigr)\nonumber\\[-8pt]\\[-8pt]
&=&e^{f(t)p -g(t)p^2}.\nonumber
\end{eqnarray}
Assume now that $t\le0$, and for $\kappa>1$, $p<\frac{\kappa
-1}{\kappa}$.
We claim that for \mbox{$0\leq x\le\frac{\kappa-1}{\kappa}$},
%
%
\begin{equation}
\label{folk-3} \exp\biggl( {-x-\frac{\kappa}{2}x^2} \biggr)\le1-x.
\end{equation}
Indeed, we have an additional inequality (iii)
$1-x+\frac{x^2}{2}\ge\exp(-x)$ when $x\ge0$. Note also that
\[
\biggl(1+\frac{\kappa}{2}x^2 \biggr) (1-x)\ge1-x+
\frac
{x^2}{2}\quad\Longleftrightarrow\quad x\le\frac{\kappa-1}{\kappa}.\vadjust{\goodbreak}
\]
Thus
\[
e^{{\kappa}x^2/{2}}(1-x)\ge\biggl(1+\frac{\kappa}{2}x^2 \biggr)
(1-x) \ge1-x+\frac{x^2}{2}\ge e^{-x}.
\]
Now, set $x=-p(e^t-1)\ge0$, so that $x \le\frac{\kappa-1}{\kappa}$.
We obtain
%
%
\begin{eqnarray}
\label{folk-5} E \bigl[ {\exp\bigl( {t\bigl(X-E[X]\bigr)} \bigr)} \bigr
]&\ge&\exp
\biggl( {-tp+p\bigl(e^t-1\bigr)- \frac{\kappa}{2} p^2
\bigl(e^t-1\bigr)^2} \biggr) \nonumber\\[-8pt]\\[-8pt]
&=& e^{f(t)p-{\kappa}g(t)p^2/{2}}.\nonumber
\end{eqnarray}
Inequalities (\ref{pos-folk}) and (\ref{neg-folk}) follow
(\ref{folk-2}) and (\ref{folk-5}).
\end{pf}

%

\section{The flashing process}\label{sec-flashing}
In this section, we construct the flashing process, and state
the crucial ``uniform hitting property.''
We then present a useful equivalent construction in terms of exploration
waves.
Finally, we explain the coupling of
Lemma~\ref{theo-coupling}, but postpone its proof to the \hyperref
[app]{Appendix}.

\subsection{Construction of the process}\label{buildbuild}

\subsubsection*{Partitioning the lattice}
We are given a sequence $\H=\{h_n, n\in\N\}$.
We partition the lattice into shells $({\cal S}_j \dvtx j \geq0)$.
For an illustration, see Figure~\ref{figflashing}.
For a given parameter $h_0 > 0$,
the first shell ${\cal S}_0$ is the ball ${\mathbb B}(0,h_0)$.
For $j\ge1$, shell j is the annulus
[see its definition (\ref{annulus-dfn})]
\[
{\cal S}_j = {\mathbb A}(r_j - h_j,
r_j + h_j),
\]
where $\{r_j, j\ge1\}$ is defined inductively by
$r_1 = h_0 + h_1$, and for $j\ge1$,
\[
r_{j+1} - h_{j+1} = r_j + h_j.
\]
In Section~\ref{sec-results}, we need that (o) $\H$ is increasing,
(i) $j\mapsto h_j/r_j$ is
decreasing and (ii) $h_j=O(r_j^{1/3})$. These properties
are a straightforward consequence of our hypothesis $h_j\le h_{j+1}\le
(1+\frac{1}{2j})h_j$. Actually
we will only need these properties,
and our hypothesis is no more than a sufficient condition.

We also define
\[
\Sigma_0 =\{0\} \quad\mbox{and}\quad\Sigma_j = \partial{\mathbb
B}(0,r_j),\qquad j\geq1.
\]

\subsubsection*{Flashing times}
The key feature we expect from
the \textit{flashing process} is that its covering
mechanism be simple. More precisely, our construction is guided by
property (ii) of the \hyperref[intro]{Introduction} which states that
\textit{the flashing position}, \textit{in a shell}, \textit{is
essentially uniform
over the cell an explorer first hits upon entering the shell}.
Thus, we need to define together \textit{cells} and
\textit{flashing times} to realize
property (ii). It is important that all sites of a shell can be chosen
as \textit{flashing sites} with about the same frequency.
In this respect, let us remark that a cell\vadjust{\goodbreak} in shell ${\cal S}_j$
cannot be a ball of radius $h_j$ centered on $\Sigma_j$. Indeed,
if this were the case,
sites at a distance about $h_j$ would be in much fewer cells
than sites of $\Sigma_j$, and this would fail to make the covering
of a shell uniform. We find it convenient to build a cell with
a \textit{mixture} of balls and annuli. A (random) flag $Y_j$ tells
the explorers whether it flashes upon exiting either a sphere or
the boundary of an annulus, whose distance from $\Sigma_j$ is
governed with a random radius $R_j$ of appropriate density. Also,
to allow for the possibility of flashing on its hitting position
on $\Sigma_j$, we introduce an additional flag~$X_j$.

More precisely,
consider $\{X_j,Y_j,j\ge0\}$ a sequence of independent Bernoulli
variables such that
\[
P(X_j = 1) = 1 - P(X_j = 0) = \frac{1}{h_j^d}
\]
and
\[
P(Y_j = 1) = 1 - P(Y_j = 0) = \cases{ %
1, &\quad
if $j = 0$,
\vspace*{2pt}\cr
\frac{1}{2}, &\quad if $j\geq1$.}
\]
Consider also a sequence of continuous independent variables
$\{R_j,j\ge0\}$ each of which has density $g_j\dvtx[0,h_j]\to\R^+$ with
%
%
\begin{equation}
\label{density} g_j(h)=\frac{dh^{d-1}}{h_j^d}.
\end{equation}
For $j\ge0$, and $z_j$ in $\Sigma_j$,
let $S$ be a random walk starting in $z_j$, an define
a stopping time $\sigma$ as follows.
If $R_j = h$ for some $h \leq h_j$, then
\[
\sigma= \cases{ %
0, &\quad if $X_j = 1$,
\vspace*{1pt}\cr
H\bigl({\mathbb
B}\bigl(z_j, h\wedge\bigl(r_j + h_j -
\|z_j\|\bigr)\bigr)^c\bigr), &\quad if $X_j = 0$ and
$Y_j = 1$,
\vspace*{2pt}\cr
H\bigl({\mathbb A}(r_j - h,
r_j + h)^c\bigr), &\quad if $X_j = 0$ and
$Y_j = 0$.}
\]
We set $H_j = H(\Sigma_j)$, and we define the stopping times
$(\sigma_j \dvtx j\geq0)$ as
\[
\sigma_j = H_j + \sigma(S\circ\theta_{H_j}),
\]
where $\theta$ stands for the usual time-shift operator.
For $j \geq0$ we note that, by construction, $S(t)\in{\cal S}_j$
for all $t$ such that $H_j \leq t < \sigma_j$ and we say that
$\sigma_j$ is a \textit{flashing time} when $S(\sigma_j)$
is contained in the intersection between ${\cal S}_j$
and the cone with base $B(S(H_j), h_j / 2)$.
We call such an intersection \textit{a cell centered at} $S(H_j)$,
that we denote $\C(S(H_j))$.
In other words, for any $z\in\Sigma_j$
%
%
\begin{equation}
\label{def-cell} {\cal C}(z) = {\cal S}_j \cap\bigl\{ x \in{\mathbb
R}^d \dvtx\exists\lambda\geq0, \exists y\in B(z,
h_j/2), x = \lambda y \bigr\}.
\end{equation}

\subsubsection*{The uniform hitting property}
The main property of the hitting time $\sigma$ constructed above
is the following proposition, which yields property (ii) of
the flashing process to be defined soon.\vadjust{\goodbreak}
%
%
\begin{itproposition}\label{prop-uniform}
There are two positive constants $\alpha_1 < \alpha_2$,
such that, for $h_0$ large enough, $j\geq0$, $z_j\in\Sigma_j$, and
$z^*\in\C(z_j)$.
%
%
\begin{equation}
\label{uniform-inequality} \frac{\alpha_1}{h_j^d} \leq\P_{z_j} \bigl( S(
\sigma)=z^* \bigr) \leq\frac{\alpha_2}{h_j^d}.
\end{equation}
\end{itproposition}
The proof of Proposition~\ref{prop-uniform} is given in
Section~\ref{sec-harmonic}.

\subsubsection*{The flashing process}
Consider a family of $N$ independent random walks
$(S_i\dvtx1\leq i \leq N)$
with their stopping times $(H_{i,j}, \sigma_{i,j} \dvtx j\geq0)$.
Let also $z_{i,j}=S_i(H_{i,j})$ be the first hitting position
of $S_i$ on $\Sigma_j$.

We define the cluster inductively. Set
$A^*(0) = \varnothing$. For $i\geq1$,
we define $\tau_i^*$ as the first flashing time
associated with $S_i$ when the
explorer stands outside $A^*(i -1)$. In other words,
\[
\tau_i^* =\min\bigl\{ \sigma_{i,j}\dvtx j \geq0,
S_i(\sigma_{i,j}) \in{\cal C}(z_{i,j})\cap
A^*(i-1)^c \bigr\}
\]
and
\[
A^*(i) = A^*(i-1)\cup\bigl\{S_i\bigl(\tau_i^*\bigr)
\bigr\}.
\]

\subsection{Exploration waves}\label{sec-wave}

Rather than building $A^*(N)$
following the whole journey of one explorer
after another, we can build $A^*(N)$
as an increasing union of clusters formed
by stopping explorers on successive shells. 
Similar wave constructions are introduced
in~\cite{lawler92} and~\cite{lawler95}.
We use this alternative construction
in the proof of Propositions~\ref{theo-flashing}
and~\ref{theo-optimal}.

We denote by $\xi_k\in({\mathbb Z}^d)^N$
the explorers positions after the $k$th wave.
We denote by ${\cal A}^*_k(N)$
the set of sites where settled explorers are
after the $k$th wave. Our inductive construction will be such that
\[
\xi_k(i)\notin\Sigma_k \quad\Leftrightarrow\quad
\xi_k(i)\in\bigcup_{j<k}{\cal S}_j \quad\Leftrightarrow\quad
\xi_k(i)\in{\cal A}^*_k(N).
\]
For $k = 0$ we set $\xi_0(i)= 0$, and
${\cal A}^*_0(i) = \varnothing$, for $1\leq i\leq N$.
Assume that for $k\geq0$, ${\cal A}^*_{k}(i)$ is built for $i=0,\ldots,N$.
We set ${\cal A}_{k+1}^*(0) = {\cal A}^*_{k}(N)$.
For $i$ in $\{1,\ldots,N\}$, we set the following:
\begin{itemize}
\item If $\xi_k(i)\notin\Sigma_k$, then
\[
\xi_{k+1}(i) = \xi_{k}(i) \in\bigcup_{j<k}{\cal
S}_j\quad\mbox{and}\quad{\cal A}^*_{k+1}(i) = {\cal
A}^*_{k+1}(i-1).
\]
\item
If $\xi_k(i)\in\Sigma_k$
and $S_i(\sigma_{i,k})\in{\cal C}(z_{i,k})\cap{\cal A}^*_k(i-1)^c$, then
\[
\xi_{k+1}(i) = S_i(\sigma_{i,k})\in{\cal
S}_k \quad\mbox{and}\quad{\cal A}^*_{k+1}(i) = {\cal
A}^*_{k+1}(i-1) \cup\bigl\{S_i(\sigma_{i,k})
\bigr\}.
\]
\item
If $\xi_k(i)\in\Sigma_k$
and $S_i(\sigma_{i,k})\notin{\cal C}(z_{i,k})\cap{\cal
A}^*_k(i-1)^c$, then
\[
\xi_{k+1}(i) = S_i(H_{i,k+1})\in
\Sigma_{k+1}\quad\mbox{and}\quad{\cal A}^*_{k+1}(i) = {\cal
A}^*_{k+1}(i-1).
\]
\end{itemize}
In words, for each $k\geq1$, during the $k$th wave of exploration,
the unsettled explorers move one after the other
in the order of their labels until either settling in ${\cal S}_{k-1}$,
or reaching $\Sigma_{k}$ where they stop.
We then define ${\cal A}^*(N)$ by
\[
{\cal A^*}(N) = \bigcup_{k\geq1}{\cal
A}^*_k(N).
\]
We explain now why this construction yields the same
cluster as our previous definition.
An explorer cannot settle inside a shell it has left,
and thus cannot settle in any shell ${\cal S}_j$
with $j<k$ if it reaches ${\Sigma}_k$. Now, since
each wave of exploration is organized according to the label ordering,
the fact that an explorer has to wait for the following explorers
before proceeding its journey beyond $\Sigma_k$ does not interfere
with the site where it eventually settles.

\subsection{Coupling internal DLA and flashing processes}
\label{sec-coupling}

\subsubsection*{Proof of Lemma \protect\ref{theo-coupling}}
For each positive integer $N$, we build a coupling between $A(N)$ and
$A^*(N)$.
We first describe the main features of our coupling in words.
Its precise definition is postponed to the \hyperref[app]{Appendix}.

We launch $N$ independent random walks, and build inductively the
associated clusters $A(1)$, $A(2),\ldots,A(N)$. In doing so, we use the
increments of these random walks to define, step by step, $N$ flashing
trajectories $S^*_1,\ldots,S^*_N$ up to some times $\bar
t_1,\ldots,\bar t_N$. Let us describe informally step $i+1$ of the
induction. Assume that $S_1^*,\ldots,S_i^*$ are defined up to some
times $t_1\leq\bar t_1,\ldots,t_i\leq\bar t_i$, and that each site of
$A(i)$ is covered by exactly one $S_k^*(t_k)$ with $1\leq k\leq i$. We
can think of $S_1^*(t_1),\ldots,S_i^*(t_i)$ as the positions of stopped
flashing explorers, some of them stopped at one of their flashing
times---say on \textit{blue} sites---some of them not---say on \textit{red}
sites. Then, we add the $i+1$th explorer and flashing explorer. We set
$S_{i+1}^*(0) = S_{i+1}(0) = 0$. We add new increments both to
$S_{i+1}$ and to the trajectory of one flashing explorer, say with
label $j$ in $\{1;\ldots;i+1\}$, in such a way that the current
position of the walker $i+1$ and that of the flashing explorer $j$
coincide. The label $j$ is defined inductively as follows. Initially,
$j = i+1$. Assume now that the walker $i+1$ flashes on a red or blue
site inside $A(i)$. This site is occupied by exactly two stopped
flashing explorers, $j$ and $j'$ [and all other red and blue sites of
$A(i)$ are occupied by exactly one flashing explorer]. Since flashing
explorers can settle at their flashing times, it makes sense, when $j$
is flashing, to add the next increment to the trajectory of flashing
explorer $j'$ rather than $j$. We do so in two cases, first, when this
happens on a red site. In this case, we turn blue that site since $j$
is stopped at a flashing time. Second, when this happens on a blue
site, say $z$, \textit{and} $j' > j$. Note that in this case, both
explorers flash on $z$, but explorer $j$ reaches $z$ before explorer
$j'$ when launched in their label order. Our choice is such that the
eventual cluster $A^*(N)$ has the correct law.\vadjust{\goodbreak} In all other cases, we
keep adding the increments of $S_{i+1}$ to the same flashing
trajectory. It is important to note that the value of the increment
does not depend on the index of the trajectory we choose to extend.
Walker $i+1$ eventually steps outside $A(i)$, say on $z^*$, while
following a flashing trajectory, say the $j$th one. We stop the $j$th
flashing trajectory on $z^*$, and paint $z^*$ blue or red according to
whether $z^*$ is one of its flashing sites or not.

When the last walker steps outside $A(N-1)$, we have
%
%
\begin{equation}
\label{afterloop1} A(N) = \bigl\{S_1^*(
\bar t_1);\ldots;S_N^*(\bar t_N)\bigr\}
\qquad\mbox{with } \bigl|A(N)\bigr| = N.
\end{equation}
To define $A^*(N)$ we launch again, in their label's order,
the flashing explorers from their current positions
(possibly some or none of them since some or all of
them can already have reached their settling position).
We then get
%
%
\begin{eqnarray}
\label{afterloop2} A^*(N) = \bigl\{S_1^*
\bigl(\tau_1^*\bigr);\ldots;S_N^*\bigl(
\tau_N^*\bigr)\bigr\} \nonumber\\[-8pt]\\[-8pt]
&&\eqntext{\mbox{with } \bigl|A^*(N)\bigr| = N \mbox{ and }
\tau_k^*\geq\bar t_k\qquad\mbox{for all $k$}.}
\end{eqnarray}

\subsubsection*{Proof of Corollary \protect\ref{cor-coupling}}
Since a flashing explorer that visited
some site beyond a given shell cannot settle
in that shell, the one-to-one map
%
%
\begin{equation}
\label{dfnpsi} \psi_N\dvtx S_k^*(\bar
t_k) \in A(N)\mapsto S_k^*\bigl(\tau_k^*
\bigr) \in A^*(N),\qquad k=1,\ldots,N,
\end{equation}
satisfies, for all $k$ and $l$,
%
%
\begin{equation}
\label{desired} S_k^*(\bar t_k) \notin\bigcup
_{m<l}{\cal S}_m \quad\Rightarrow\quad S_k^*
\bigl(\tau_k^*\bigr) = \psi_N\bigl(S_k^*(\bar
t_k)\bigr) \notin\bigcup_{m<l}{\cal
S}_m.
\end{equation}

Thus, for all $N\geq0$ there is a coupling and a one-to-one map $\psi_N$
between $A(N)$ and $A^*(N)$ such that for all $k\geq1$,
%
%
\begin{equation}
\label{inclus-1} \psi_N \bigl(A(N)\cap{\mathbb B}^c_{r_k+h_k}
\bigr) \subset A^*(N)\cap{\mathbb B}^c_{r_k+h_k}.
\end{equation}
Inclusion (\ref{inclus-1}) has two important consequences:
\begin{longlist}[(a)]
\item[(a)] If $A^*(N)\subset{\mathbb B}_{r_k+h_k}$,
then $A(N)\subset{\mathbb B}_{r_k+h_k}$. Indeed,
any site in $A(N)$ outside
${\mathbb B}_{r_k+h_k}$ produces, through $\psi_N$,
a site in $A^*(N)$ outside ${\mathbb B}_{r_k+h_k}$.
\item[(b)] If ${\mathbb B}_{r_k+h_k}\subset A^*(N)$, then
${\mathbb B}_{r_k+h_k}\subset A(N)$.
Indeed, those sites in $A(N)$ that are mapped through $\psi_N$
on $A^*(N)\cap\B_{r_k+h_k} = \B_{r_k+h_k}$
are necessarily contained in ${\B}_{r_k+h_k}$.
Since their number is $|{\mathbb B}_{r_k+h_k}|$ and
$\psi_N$ is one-to-one,
they completely cover ${\mathbb B}_{r_k+h_k}$.
\end{longlist}

\section{Fluctuations}\label{sec-results}
In this section, we prove Propositions~\ref{theo-flashing}
and~\ref{theo-optimal}.
To do so we use the construction in terms of exploration waves
of Section~\ref{sec-wave}. Thus, we think of the growing cluster
as evolving in discrete time,
where time counts the number of exploration waves. The proofs
in this section rely on potential theory estimates which we
have gathered in Section~\ref{sec-harmonic}, for the ease of reading.

\subsection{Tiles}
We recall that we have defined a \textit{cell} of $\S_j$
in (\ref{def-cell}), as the intersection of a cone with $\S_j$.
We need also a smaller shape.\vadjust{\goodbreak}
We define, for any $z_j$ in $\Sigma_j$, and for a small $\varepsilon_0$ to
be defined later,
%
%
\begin{equation}\label{dfn-cell2}
\tilde{\cal C}(z_j) = {\cal S}_j \cap\bigl\{ x \in{
\mathbb R}^d \dvtx\exists\lambda\geq0, \exists y\in
B(z_j, \varepsilon_0h_j), x = \lambda y \bigr
\}.
\end{equation}
As in Lemma 12 in~\cite{lawler95}, concerning locally finite coverings,
we claim that, for $h_0$ large enough,
there exist a positive constant $K_F$,
and, for each $j\geq0$, a subset $\tilde\Sigma_j$
of $\Sigma_j$ such that
%
%
\begin{equation}
\label{local-covering} \forall y\in\S_j\qquad \bigl\llvert\bigl\{ {z\in
\tilde\Sigma_j\dvtx y\in\tilde\C(z)} \bigr\}\bigr\rrvert\le
K_F \quad\mbox{and}\quad\S_j = \bigcup
_{z_j\in\tilde\Sigma_j} \tilde{\cal C}(z_j).
\end{equation}
For any $z_j\in\tilde\Sigma_j$, we call \textit{tile centered at} $z_j$,
the intersections of $\tilde{\cal C}(z_j)$ with $\Sigma_j$.
We denote by $\T(z_j)$ a tile centered at $z_j$, and
by ${\cal T}_j$
the set of tiles associated with the shell ${\cal S}_j$.
%
%
\begin{equation}
{\cal T}_j = \bigl\{ {\cal T}(z_j)\dvtx
z_j\in\tilde\Sigma_j \bigr\}.
\end{equation}
We choose $\varepsilon_0$ to satisfy two properties.
First,
for any $z\in\S_j$, there is $\tilde z_j\in\tilde\Sigma_j$ such that
%
%
\begin{equation}
\label{feature-coupon} z\in\bigcap_{y\in{\cal T}(\tilde z_j)}\C(y).
\end{equation}
This is ensured by the choice of a small enough $\varepsilon_0$.
Indeed, let $z_j\in\Sigma_j$ be a site realizing
the minimum of $\{\|z-y\|\dvtx y\in\Sigma_j\}$. There is
$\lambda>0$ and $u\in B(z_j,1)$, such that $z=\lambda u$.
Now, there is $\tilde z_j\in\tilde\Sigma_j$ such that
$\|\tilde z_j-z_j\|<\varepsilon_0h_j$, and for any $y\in\T(\tilde z_j)$,
we have $\|y-z_j\|<2\varepsilon_0 h_j$. Thus, for
$\varepsilon_0$ small enough so that $1+2\varepsilon_0 h_j\le h_j/2$,
\[
\forall y\in\T(\tilde z_j)\qquad\|u-y\|\le\|u-z_j\|+
\|z_j-y\|\le1+2\varepsilon_0 h_j\le
\frac{h_j}{2},
\]
which implies (\ref{feature-coupon}).
Second, the size of a tile should be such that for some $\kappa>1$,
for any $j\ge1$, and any tile $\T\in\T_j$
%
%
\begin{equation}
\label{cond-kappa} \sup_{z\in\B(0,r_j-h_j)} \P_z \bigl( {S\bigl(H(
\Sigma_j)\bigr)\in\T} \bigr)\le\frac{\kappa-1}{\kappa}.
\end{equation}
Inequality (\ref{cond-kappa}) follows from Lemma 5(b) of
\cite{lawler92} (or Lemma~\ref{lem-alex} below)
which for a constant $J_d$ yields
\[
\sup_{z\in\B(0,r_j-h_j)} \P_z \bigl( {S\bigl(H(\Sigma_j)
\bigr)\in\T} \bigr)\le J_d\frac{|\T|}{h_j^{d-1}}.
\]
The choice of $\varepsilon_0$ is such that $J_d |\T|\le
\frac{\kappa-1}{\kappa} h_j^{d-1}$.

\subsection{Bounding inner fluctuations}

For $n\geq0$, we take $N = |{\mathbb B}_n|$, we recall that
${\cal A}^*_k(N)\subset{\cal A}^*_{k+1}(N)$ for $k\in\N$, and
${\cal A}^*(N) = \bigcup_{k\geq1} {\cal A}^*_k(N)$.
We consider
%
%
\begin{equation}
T^* = \min\biggl\{ k\geq1 \dvtx\bigcup_{j<k}{\cal
S}_j \not\subset{\cal A}^*_k(N) \biggr\}.\vadjust{\goodbreak}
\end{equation}
Note that ${\cal A}^*_k(N)\subset\bigcup_{j<k}{\cal S}_j$, so that
$T^*$ is the first time $k$
when the $k$th wave does not cover all its allowed space.
We recall that time counts the number of exploration waves.

For the flashing process if $\bigcup_{j<k}\S_j
\not\subset{\cal A}^*_k(N)$, then for any $k'>k$, we have
$\bigcup_{j<k}\S_j \not\subset{\cal A}^*_{k'}(N)$, so that $T^*$
is also the shell label where the first hole of ${\cal A}^*(N)$
appears. We have, for $l$ with $r_l <n$,
%
%
\begin{equation} \label{est0}
P\bigl(T^*\le l\bigr)= P \bigl({\mathbb B}(0,r_l +
h_l)\not\subset{\cal A}^* (N)\bigr) \leq\sum
_{k\leq l} P\bigl(T^* = k+1\bigr).
\end{equation}
In this section, we estimate from above the probability $P(T^* = k+1)$
assuming $r_k < n$.

For $k\geq1$ and $\Lambda\subset\Sigma_k$,
we call $W_k(\Lambda)$ the number of unsettled explorers
that stand in $\Lambda$ after the $k$th wave,
that is,
%
%
\begin{equation}
W_k(\Lambda) = \sum_{i = 1}^N
\id_{\Lambda} \bigl(\xi_k(i) \bigr).
\end{equation}
We now look at the
\textit{crossings} of tiles of $\T_k$. On the one hand, we will use that
if $W_k({\cal T})$ is \textit{large}, then it is unlikely that a hole
appears in the cell containing $\T$ during the $k+1$th-wave.
We use for this purpose the fact that covering for the flashing
process is similar to filling an album for a coupon-collector model.
On the other hand, if $r_k$ is
\textit{small}, it is unlikely that $W_k({\cal T})$ is \textit{small}.
We now make precise what we intend by \textit{small} and \textit{large}.
For any positive constant $\xi$, we write
%
%
\begin{eqnarray}
\label{est1} P \bigl(T^* = k+1 \bigr) &=& P \bigl( T^* = k + 1, \forall
{\cal T}
\in{\cal T}_k, W_k({\cal T}) \geq\xi\bigr)\nonumber\\
&&{} + P \bigl( T^*
= k + 1, \exists{\cal T}\in{\cal T}_k, W_k({\cal T}) <
\xi\bigr)
\nonumber\\[-8pt]\\[-8pt]
&\le& P \bigl( {T^*=k+1 | \forall{\cal T} \in{\cal T}_k,
W_k({\cal T}) \geq\xi} \bigr)\nonumber\\
&&{}+ P \bigl(\exists{\cal T}\in{\cal
T}_k, W_k({\cal T}) < \xi\bigr).\nonumber
\end{eqnarray}

\subsubsection*{A coupon-collector estimate}
The first term in the right-hand side of (\ref{est1})
is bounded using a simple coupon-collector argument.
Indeed, the event $\{T^* = k + 1\}$ implies that there is an
uncovered site in $\S_k$, say $z$, when explorers stopped in $\Sigma_k$
are released.
By (\ref{feature-coupon}), there is $z_k\in
\tilde\Sigma_k$, such that $z$ is a possible settling
position of all explorers stopped in $\T(z_k)$. Now, knowing that
$\{W_k(\T(z_k)) \geq\xi\}$, Proposition~\ref{prop-uniform}
tells us that the probability of not covering this site is
less than $(1-\alpha_1/h_k^d)$ to the power $\xi$.
In other words,
\[
P \bigl( T^* = k + 1 | \forall{\cal T}\in{\cal T}_k,
W_k({\cal T}) \geq\xi\bigr) \leq|{\cal S}_k| \biggl( 1
- \frac{\alpha_1}{h_k^d} \biggr)^{\xi} \leq|{\cal S}_k| \exp
\biggl( {-\alpha_1 \frac{\xi}{h_k^d}} \biggr).
\]
Henceforth, we set
%
%
\begin{equation}
\label{def-h} \xi=Ah^d\log(n) \qquad\mbox{with } h=\sup\{h_k
\dvtx r_k\le n\}
\end{equation}
and $A$ large enough so that
%
%
\begin{eqnarray}
\label{inter-9}
&&
\sum_{k\dvtx r_k<n} P \bigl( T^* = k + 1 |
\forall{\cal T}\in{\cal T}_k, W_k({\cal T})
\geq\xi\bigr)\nonumber\\[-8pt]\\[-8pt]
&&\qquad\le|{\mathbb B_n}| \exp( {-\alpha_1 A\log n } )\le
\frac{1}{n^2}.\nonumber
\end{eqnarray}

\subsubsection*{Estimating $\{W_k({\cal T}) < \xi\}$}
For any $\T\in\T_k$, we consider
the counting variable $L_k(\T)=M(\B(0, r_k - h_k), r_k, \T)$,
and define
%
%
\begin{eqnarray}
\label{identity}
&&M_k(\T)=W_k(\T) + M\bigl(A^*_k,r_k,
\T\bigr)\nonumber\\[-8pt]\\[-8pt]
&&\eqntext{\mbox{so that } M_k(\T)\stackrel{\mathrm{law}} {=}M(N
\id_{\{0\}}, r_k, \T).}
\end{eqnarray}
The idea of defining $M_k$ and $L_k$ (for the internal DLA process),
and bounding $W_k$ by $M_k-L_k$, is introduced in~\cite{lawler92}.
Our main observation is that $L_k(\T)$ is independent of $W_k(\T)$,
and
\[
W_k(\T) + L_k(\T)\ge M_k(\T).
\]
As a consequence, for any positive constants $t$ and $\xi$
(and with the notation $\bar X=X-E[X]$),
\begin{eqnarray*}
P \bigl( {W_k(\T)<\xi} \bigr) &\le& e^{t \xi} \times E \bigl[ {
\exp\bigl( {-t W_k(\T)} \bigr)} \bigr] =e^{t \xi}
\frac{E [ {\exp( {-t ( {W_k(\T)+L_k(\T
)} )} )} ]} {
E [ {\exp( {-t L_k(\T)} )} ]}
\\
&\le& \exp\bigl( {-t \bigl( {E \bigl[ {M_k(\T)-L_k(\T)}
\bigr]-\xi} \bigr)} \bigr) \times\frac{E [ {\exp( {-t ( {\bar M_k(\T)}
)} )} ]} {
E [ {\exp( {-t \bar L_k(\T)} )} ]}.
\end{eqnarray*}
Using Lemma~\ref{lem-folk} with condition (\ref{cond-kappa}),
we obtain
\begin{eqnarray*}
\log P \bigl( {W_k(\T)<\xi} \bigr)
&\le& -t \bigl( {E \bigl[
{M_k(\T)-L_k(\T)} \bigr]-\xi} \bigr)+ f(-t)E \bigl[
{M_k(\T)-L_k(\T)} \bigr]
\\
&&{} +\frac{\kappa}{2}g(-t) \sum_{y\in\B(0, r_k - h_k)}P_y^2
\bigl( {S\bigl(H(\Sigma_k)\bigr)\in\T} \bigr).
\end{eqnarray*}
We now proceed in two steps.
We show in step 1 that for some constant $\kappa'$,
\[
E \bigl[ {M_k(\T)-L_k(\T)} \bigr]\ge
\kappa' \bigl(n^d-(r_k-h_k)^d
\bigr) \frac{h_k^{d-1}}{r_k^{d-1}}.
\]
Since $\{h_k / r_k, k \geq0\}$ is nonincreasing,
it follows that there is a constant $\kappa_1>0$
such that, for all $\alpha>0$, and
$k_\alpha:= \sup\{j\in{\mathbb N}\dvtx r_j < n-\alpha h \log n\}$,
where $h$ is defined in~(\ref{def-h}), we have
\[
\inf_{k\leq k_\alpha} E \bigl[ {M_k(\T)-L_k(\T)} \bigr]
\ge\kappa' \bigl(n^d-(n-h)^d\bigr)
\frac{h^{d-1}}{n^{d-1}}\ge\kappa_1 \alpha h^d\log n.
\]
Now, if we choose $\xi$ as in (\ref{def-h}), with $\alpha= 2A/\kappa_1$
and $k^*= k_\alpha$, that is,
\[
k^*:=\sup\biggl\{ {j\in{\mathbb N}\dvtx r_j\le n-\frac{2A}{\kappa_1}h
\log(n)} \biggr\},
\]
then, we get, for all $k\leq k^*$,
%
%
\begin{equation}
\label{mean-xi} E \bigl[ {M_k(\T)-L_k(\T)} \bigr]\ge2\xi.
\end{equation}

We show in step 2, that for a constant $C$
depending on the dimension only
%
%
\begin{equation}
\label{log-4} \sum_{y\in\B(0, r_k - h_k)}\P_y^2
\bigl( {S\bigl(H(\Sigma_k)\bigr)\in\T} \bigr)\le C E \bigl[
{M_k(\T)-L_k(\T)} \bigr].
\end{equation}
Suppose for a moment that steps 1 and 2 hold. Since,
for some $c> 0$,
$\max(f(-t)$, $g(-t))\le c t^2$ when $t\le1$, there is $c'>0$ such that
for $k\le k^*$
%
%
\begin{eqnarray}
\label{log-main}
&&
\log P \bigl( {W_k(\T)<Ah^d\log(n)}
\bigr)\nonumber\\
&&\qquad\le\inf_{0\le t\le1} \biggl( {-{t}+c\biggl(1+\frac{C\kappa}{2}
\biggr)t^2} \biggr)E \bigl[ {M_k(\T)-L_k(\T)}
\bigr]
\\
&&\qquad\le-c' E \bigl[ {M_k(\T)-L_k(\T)} \bigr]
\le-2c' Ah^d\log(n).\nonumber
\end{eqnarray}
Now, using (\ref{est1}), (\ref{inter-9}) and (\ref{log-main})
for $A$ large enough, we have
\[
\sum_{k< k^*} P\bigl(T^* = k \bigr)\le
\frac{2}{n^2}.
\]
Borel--Cantelli's lemma yields then
the inner control of Proposition~\ref{theo-flashing}.

\subsubsection*{Step 1}
We invoke Corollary~\ref{cor-gauss},
with $n=r_k$, and $\Delta_n=h_k$
[the hypotheses $h_k=O(r_k^{1/3})$ and $h_k$ large enough
hold here, as seen in the first paragraph of Section~\ref{buildbuild}].
We have
for some positive constants $\kappa'$, $K$ and for $n$ large enough,
%
%
\begin{eqnarray}
\label{condtalt}\quad E \bigl[ M_k({\cal T}) - L_k({\cal T})
\bigr] & = & E\bigl[M\bigl(\bigl(|{\mathbb B}_n|-|{\mathbb
B}_{r_k-h_k}|\bigr)\id_{0},r_k,\T\bigr)\bigr]
\nonumber
\\
&&{} + E\bigl[M\bigl(|{\mathbb B}_{r_k-h_k}|\id_{0},r_k,
\T\bigr)\bigr]- E\bigl[M({\mathbb B}_{r_k-h_k},r_k,\T)\bigr]
\nonumber
\\
&\ge& \bigl( |{\mathbb B}_n| - |{\mathbb B}_{r_k-h_k}| \bigr) {\mathbb
P}_0 \bigl( S(H_k)\in{\T} \bigr) -K h_k^{d-1}
\\
&\geq& 2\kappa'\bigl(n^d - (r_k-h_k)^d
\bigr)\frac
{h_k^{d-1}}{r_k^{d-1}}-Kh_k^{d-1}
\nonumber
\\
&\geq& \kappa'\bigl(n^d - (r_k-h_k)^d
\bigr)\frac{h_k^{d-1}}{r_k^{d-1}}\nonumber
\end{eqnarray}
for $r_k\leq n$ and $h_0$ large enough.

\subsubsection*{Step 2}
By Lemma~\ref{lem-alex} below,\vspace*{1pt} there is a constant $\kappa
_G$ such that,
for $y\in\B(0, r_k - h_k)$, and $z\in\tilde\Sigma_k$
\[
\P_y \bigl( {S\bigl(H(\Sigma_k)\bigr)\in\T(z)} \bigr)
\le\frac{\kappa_G
|\T(z)|} {
\|z-y\|^{d-1}}.
\]
Therefore
%
%
\begin{equation}
\label{log-6} \sum_{y\in\B(0, r_k - h_k)}\P_y^2
\bigl( {S\bigl(H(\Sigma_k)\bigr)\in\T(z)} \bigr)\le\sum
_{j\dvtx h_k\le j\le2r_k} \sum_{y\dvtx j\le|z-y|<j+1}
\frac{\kappa_G^2 |\T(z)|^2}{j^{2(d-1)}}.\hspace*{-35pt}
\end{equation}
For a constant $C_d$, we bound $|\{y\dvtx k\le|z-y|<k+1\}|\le C_d k^{d-1}$.
Thus,
%
%
\begin{eqnarray}
\label{log-7}
&&
\sum_{y\in\B(0, r_k - h_k)}\P_y^2
\bigl( {S\bigl(H(\Sigma_k)\bigr)\in\T(z)} \bigr)\nonumber\\[-3pt]
&&\qquad\le\sum
_{j\dvtx h_k\le j\le2r_k} \frac{C_d\kappa_G^2 |\T
(z)|^2}{j^{d-1}}
\\[-3pt]
&&\qquad\le C' \bigl|\T(z)\bigr|^2 \biggl( {\ind_{d=2}
\log(n)+\ind_{d>2} \frac
{1}{h_k^{d-2}}} \biggr).\nonumber
\end{eqnarray}
Since $|\T(z)|$ is of order $h_k^{d-1}$, (\ref{log-4}) holds.

\subsection{Bounding outer fluctuations}\label{sec-outer}

This section follows~\cite{lawler95} closely.
The features of the flashing
process allow for some simplification.
We keep the notation of the previous subsection.
There, we proved that for some integer~$k^*$,
which depends on $n$,
\[
P \bigl( {T^*>k^*} \bigr) =1 - \varepsilon(n)\qquad\mbox{with } \sum
_{n\geq1} \varepsilon(n) < +\infty.
\]
The integer $k^*$ is the largest such that $r_{k^*}\le n-
2A h\log(n)/\kappa_1$, for a large constant $A$
and with $h$ defined in (\ref{def-h}).
As a consequence, the following conditional law
can be seen as a slight modification of $P$:
%
%
\begin{equation}
P^*(\cdot) = P \bigl( {\cdot|T^*>k^*} \bigr).
\end{equation}

We begin by proving that under $P^*$
the probability to find some $k$
with $n\le r_k < 2n$ and some tile ${\cal T}$
in ${\cal T}_k$ with $W_k({\cal T})$
larger than or equal to $\xi' = 2A'h^d\log n$
for a large enough $A'$ decreases
faster than any given power of $n$.
First, note that on $\{T^*>k^*\}$,
%
%
\begin{equation}
W_k({\cal T}) + L_k^*({\cal T})\leq M_k({
\cal T})\qquad\mbox{with } L_k^*= M \bigl( {\mathbb B}(0,
r_{k^*}-h_{k^*}), r_k, {\cal T} \bigr).\hspace*{-30pt}
\end{equation}
Our key observation is that the pair $(W_k({\cal T}), \ind_{\{T^*>k^*\}
})$ is
independent of~$L_k^*$. Thus, for any $t>0$,
\begin{eqnarray*}
P \bigl( {W_k({\cal T}) \geq\xi', T^*>k^*} \bigr)
&\le&
e^{-t\xi'} E \bigl[ {e^{tW_k({\cal T})}\ind_{\{T^*>k^*\}}} \bigr]
\\[-3pt]
&=&
e^{-t\xi'}\frac{
E [ {\exp( {t ( {W_k({\cal T})+ L_k^*} )}
)\ind_{\{T^*>k^*\}}} ]} {
E [e^{tL_k^*} ]}
\\[-3pt]
&\le& e^{-t\xi'}\frac{E [ {e^{tM_k({\cal T})}} ]}{E
[ {e^{tL_k^*}} ]}
\\[-3pt]
&=&\exp\bigl( {-t \bigl( {\xi' - E \bigl[ {M_k({\cal
T})-L_k^*} \bigr]} \bigr)} \bigr)\times\frac{E [ {e^{t\bar M_k({\cal
T})}} ]}{E [ {e^{t\bar
L_k^*}} ]}.
\end{eqnarray*}
By Lemma~\ref{lem-folk}, we have [for $f(t)$ and $g(t)$ quadratic for $t$
small]
%
%
\begin{eqnarray}
\label{log-11}
&&
\log P \bigl( {W_k({\cal T}) \geq\xi',
T^*>k^*} \bigr) \nonumber\\[-3pt]
&&\qquad\le-t \bigl( {\xi'-E \bigl[ {M_k({
\cal T})-L_k^*} \bigr]} \bigr) +f(t) \times E \bigl[
{M_k({\cal T})-L_k^*} \bigr]
\\[-3pt]
&&\qquad\quad{}+g(t) \times\sum_{y\in\B(0,r_{k^*}-h_{k^*})} \P^2_y
\bigl( {S\bigl(H(\Sigma_k)\bigr)\in{\cal T}} \bigr).\nonumber
\end{eqnarray}
The steps are now similar to the previous proof. We first
estimate\break
$E [ {M_k({\cal T})-L_k^*} ]$. By Corollary~\ref{cor-gauss},
for some positive constant $K'$ and for $n$ large enough,
%
%
\begin{eqnarray}
\label{inner-3} E \bigl[ M_k({\cal T}) - L_k^*({\cal T})
\bigr] &\leq& K' \bigl(n^d - (r_{k^*})^d
\bigr)\frac{h_k^{d-1}}{r_k^{d-1}} + O\bigl(h_k^{d-1}\bigr)
\nonumber\\[-9pt]\\[-9pt]
&\leq& K'dn^{d-1}(n-r_{k^*})\frac{h_k^{d-1}}{r_k^{d-1}} +
O\bigl(h_k^{d-1}\bigr).\nonumber
\end{eqnarray}
Note that since $r_k\le2n$, we have
$r_k^{d-1}=o(n^{d-1}(n-r_{k^*}))$ so that $O(h_k^{d-1})$
is small compared to the first term in (\ref{inner-3}).
Since $k\mapsto h_k/r_k$ is decreasing, we have for some constant $K$
\[
E \bigl[ M_k({\cal T}) - L_k^*({\cal T}) \bigr]\le
Kh^{d}\log n.
\]
Second, we estimate the sum of $P^2_y ( {S(H(\Sigma_k))\in{\cal
T}} )$
which appears on (\ref{log-11}). We use (\ref{log-6}) again
to obtain as in (\ref{log-7}), and for a constant $C$,
\begin{eqnarray*}
&&
\sum_{y\in\B(0,r_{k^*}-h_{k^*})} \P^2_y \bigl(
{S\bigl(H(\Sigma_k)\bigr)\in{\cal T}} \bigr)\\[-3pt]
&&\qquad\le C
h_k^{2(d-1)} \biggl( {\ind_{d=2} \log(n)+
\ind_{d>2} \frac
{1}{(r_k-(r_{k^*}-h_{k^*}))^{d-2}}} \biggr).\nonumber
\end{eqnarray*}
Note that $r_k-(r_{k^*}-h_{k^*})\ge h_k$, and
since $k\mapsto h_k/r_k$ is decreasing,
we have, for $n$ large enough,
$h_k\!\le\! h_{k^*}(r_k/r_{k^*})\!\le\! h\!\times\!(2n)/(n/2) $. Thus,
for a constant~$C$,\looseness=-1
\[
\sum_{y\in\B(0,r_{k^*}-h_{k^*})} \P^2_y \bigl(
{S\bigl(H(\Sigma_k)\bigr)\in{\cal T}} \bigr)\le C \bigl( {
\ind_{d=2}h^2 \log(n)+\ind_{d>2} h^d}
\bigr)\le C h^d \log(n).
\]\looseness=0
We choose $A'=K$ to obtain, for any $t > 0$
\[
\log P \bigl( {W_k({\cal T}) \geq\xi', T^*>k^*} \bigr)
\le- \bigl( {Kt-Kf(t)-Cg(t)} \bigr) h^d \log(n).\vadjust{\goodbreak}
\]
Since we have $P(T^*>k^*)\ge1/2$ for $A$ large enough
and $K$ can be taken as large as we want, we have
that $P^*(W_k({\cal T}) \geq\xi')$ decreases faster
than any given power of $n$.

Now, let $F_k$ denote the event
that no tile ${\cal T}$ in $\Sigma_k$ contains
more than $\xi' = 2A'h^d\log n$ unsettled explorers
after the $k$th exploration wave.
We define, with the notation
of Section~\ref{sec-flashing},
$\GG_k=\sigma(\xi_0,\ldots,\xi_k)$, and note
that $F_k$ and
$\{T^*>k^*\}$ are $\GG_k$-measurable.

For any tile
${\cal T}\in\T_k$, let $z_k\in\tilde\Sigma_k$ be such that
$\T=\T(z_k)$, and denote by $\tilde{\cal C} = \tilde{\cal C}(z_k)$.
We are entitled, by Proposition~\ref{prop-uniform},
to use a coupon-collector estimate on the number of
settled explorers during the $k+1$th exploration wave.
On $F_k\cap\{T^*>k^*\}$, and
for some positive constant $K_1$,
\begin{eqnarray*}
E \bigl[\bigl|{\cal A}^*_{k+1}\cap\tilde{\cal C}\bigr| | \GG_k
\bigr] &\geq& |\tilde{\cal C}| \biggl( 1 - \biggl( 1 - \frac{\alpha
_1}{h_k^d}
\biggr)^{W_k({\cal T})} \biggr)
\\
&\geq& |\tilde{\cal C}| \biggl( 1 - \exp\biggl\{ -\alpha_1
\frac{W_k({\cal T})}{h_k^d} \biggr\} \biggr)
\\
&=& \frac{|\tilde{\cal C}|}{h_k^d} W_k({\cal T}) \frac{h_k^d}{W_k({\cal
T})} \biggl( 1 -
\exp\biggl\{ -\alpha_1\frac{W_k({\cal T})}{h_k^d} \biggr\} \biggr)
\\
&\geq& K_1 W_k({\cal T}) \inf_{x\leq2A'\log n}
\frac{1 - e^{-\alpha_1x}}{x}.
\end{eqnarray*}
We now write for some positive constant $K_2$,
\begin{eqnarray*}
\inf_{x\leq2A'\log n} \frac{1 - e^{-\alpha_1x}}{x} &\geq&
\frac
{1}{2A'\log n}
\inf_{x\leq2A'\log n} \frac{1 - e^{-\alpha_1x/2A'\log n}}{x/2A'\log n}
\\
&\geq&\frac{1}{2A'\log n} \inf_{x\leq1} \frac{1 - e^{-\alpha_1x}}{x} \ge
\frac{K_2}{\log n}.
\end{eqnarray*}
We conclude that on $F_k\cap\{T^*>k^*\}$,
%
%
\begin{equation}
\label{one-cell} E \bigl[ {\cal A}^*_{k+1}\cap\tilde{\cal C} |
\GG_k \bigr] \geq K_1K_2 \frac{W_k({\cal T})}{\log n}.
\end{equation}
Recall now that property (\ref{local-covering}) implies that $K_F |
{\cal A}^*_{k+1}\cap{\cal S}_k|\ge\sum_{z_k\in\tilde\Sigma_k} |
{\cal A}^*_{k+1}\cap\tilde\C(z_k)|$.
Thus, summing over $z_k\in\tilde\Sigma_k$
with $\C=\C(z_k)$ and $\T= \T(z_k)$
in (\ref{one-cell}),
we obtain on $F_k\cap\{T^*>k^*\}$,
\[
E \bigl[\bigl|{\cal A}^*_{k+1}\cap{\cal S}_k\bigr| |
\GG_k \bigr] \geq K\frac{W_k({\cal S}_k)}{\log n}\qquad\mbox{where } K=
\frac{K_1K_2}{K_F}.
\]
Also, since $W_k({\cal S}_k) \leq|{\mathbb B}(0,n)|$,
we have, for $n$ large enough,
\[
E \bigl[\id_{F_k\cap\{T^*>k^*\}} \bigl|{\cal A}^*_{k+1}\cap{\cal
S}_k\bigr| \bigr] \geq K\frac{E[\id_{\{T^*>k^*\}}
W_k({\cal S}_k)]}{\log n} - n^dP
\bigl(F_k^c\bigr).
\]
Since $P(T^*>k^*)\ge1/2$,
%
%
\begin{equation}
E^* \bigl[\bigl|{\cal A}^*_{k+1}\cap{\cal S}_k\bigr| \bigr] \geq K
\frac{E^*[W_k({\cal S}_k)]}{\log n} - 2 n^dP\bigl(F_k^c\bigr).
\end{equation}
In other words, noting that $|{\cal A}^*_{k+1}\cap{\cal S}_k|=
W_k({\cal S}_k)-W_{k+1}({\cal S}_{k+1})$,
%
%
\begin{equation}
\label{inner-7} E^* \bigl[ W_{k+1}({\cal S}_{k+1}) \bigr]
\leq\biggl(1 - \frac{K}{\log n} \biggr) E^* \bigl[ W_{k}({\cal
S}_{k}) \bigr] + 2n^dP\bigl(F_k^c
\bigr).
\end{equation}
By iterating (\ref{inner-7}), and using our previous estimate on
$P^*(W_k({\cal T})\geq\xi')$, we obtain that for a large enough
$\varepsilon$, $E^*[W_{l_n+ \varepsilon\log^2n}({\cal
S}_{l_n+\varepsilon\log^2 n})]$, is summable,\vspace*{2pt} when $l_n$
is the lowest index for which $r_{l_n} \geq n$. Also, the probability
(under $P$!) of seeing at least one explorer reaching the shell ${\cal
S}_{l_n+\varepsilon\log^2 n}$ is summable. Using\vspace*{1pt} the
Borel--Cantelli lemma, this yields the proof of Proposition
\ref{theo-flashing}.

\subsection{Lower bound for the deviations}\label{sec-optimal}
\subsubsection{\texorpdfstring{Proof of Proposition \protect\ref{theo-optimal}: The outer deviation}
{Proof of Proposition 1.6: The outer deviation}}
We denote by $K_n$ the largest index
such that $\S_{K_n}\subset\B(0,n)$, and by $E_n$ the event
that all explorers stopped on $\Sigma_{K_n}$, at time $K_n$,
settle afterward in one of the shells
$\{\S_j \dvtx K_n\le j< K_n+b\log(n)\}$ for some positive constant $b$,
and note that $E_n=\{A^*(N)\subset\bigcup_{j<K_n+b\log(n)}
\S_j\}$.\vspace*{1pt}

We want to find $b$ such\vspace*{1pt} that $\sum_{n\ge1} P(E_n)<\infty$.
Using that the flashing times of the different explorers
are independent, we have
\begin{eqnarray*}
P(E_n) &\le& E \biggl[P \biggl(\mbox{all explorers, stopped in $
\Sigma_{K_n}$, flash in }\bigcup_{j<K_n+b\log(n)} \S_j |
\GG_{K_n} \biggr) \biggr]
\\
&\le& E \biggl[ \biggl( \sup_{z\in\Sigma_{K_n}} P \biggl(\mbox{an explorer,
started on $z$,}\\
&&\hspace*{63pt}\mbox{flashes in }\bigcup_{j<K_n+b\log(n)} \S_j |\GG
_{K_n} \biggr)
\biggr)^{W_{K_n}(\Sigma_{K_n})} \biggr].
\end{eqnarray*}
Also, there are at least $|\S_{K_n}|$
explorers stopped on $\Sigma_{K_n}$, and
there is a positive $\varepsilon_0$ such that the probability of crossing
a given shell without flashing is larger
than~$\varepsilon_0$. Thus
\begin{eqnarray*}
P(E_n)&\le& E \Bigl[ \Bigl(1- \inf_{z\in\Sigma_{K_n}} P \bigl(\mbox{an
explorer started on $z$}\\
&&\hspace*{77.5pt} \mbox{is unsettled at time }K_n+b\log(n) \bigr)
\Bigr)^{
|\S_{K_n}|} \Bigr]
\\
&\le& \bigl( {1-\varepsilon_0^{b\log(n)}} \bigr)^{|\S_{K_n}|}.
\end{eqnarray*}
When choosing $b$ small enough, we reach $\sum_{n\ge1} P(E_n)<\infty$.

\subsubsection{\texorpdfstring{Proof of Proposition \protect\ref{theo-optimal}: The inner deviation}
{Proof of Proposition 1.6: The inner deviation}}
We recall that $K_n$ is
the largest index such that $\S_{K_n}\subset\B(0,n)$.
The rough idea here is that when we stop explorers on $\Sigma_{K_n-1}$,
there are necessarily
tiles (of $\Sigma_{K_n-1}$) containing of the order of
$h^{d-1}_{K_n-1}$ sites and which receive $h^d_{K_n-1}$ explorers.
The number of explorers on these tiles is not enough to cover
the associated cells with the \textit{coupon collector} mechanism.
We now make rigorous such an argument for shells with index
of order $K_n-\log(K_n)$.

To simplify the notation, let us first define
three positive constants $c_1,c_2$ and $c_3$ such that for any $k$ with
$n/2\le r_k\le n$, we have
%
%
\begin{eqnarray}
\label{def-constants} |\S_k|&\le& c_1 h_k
n^{d-1},\nonumber\\
\frac{h_k^d}{{\sup_{z\in\Sigma_k}}|\B(z,6h_k)\cap\Sigma_k|}
|\Sigma_k|&\ge& c_2
n^{d-1} h_k\quad\mbox{and}\\
\inf_{z\in\Sigma_k}\bigl|\tilde\C(z)\bigr|&\ge&
c_3 h_k^d.\nonumber
\end{eqnarray}
Using $\alpha_2$ given in Proposition~\ref{prop-uniform},
we define
%
%
\begin{equation}
\label{def-aA} a_n=\frac{1}{8\alpha_2} \log(h_{K_n})
\quad\mbox{and}\quad A_n= \biggl[ \frac{c_2 c_3}{4 c_1} a_n \biggr].
\end{equation}
Now we assume $h_0$ large enough to have $A_n$ a
strictly positive integer.

We wish now to consider a peel of $A_n$ shells before $\partial\B(0,n)$.
Let $I_n$ be the index of the inner shell in this peel, that is,
$r_{I_n+A_n}\le n< r_{I_n+A_n+1}$.
Since $T^*\le I_n+1$ implies that $\bigcup_{j\le I_n}\S_j\not\subset A^*(N)$,
it is enough
to show that $P(T^*>I_n+1)$ decays faster than any polynomial in $n$.

Note that the monotonicity of $k\mapsto h_k/r_k$ and $r_{I_n}\ge n/2$,
imply that $2h_{I_n}\ge h_k$ for $I_n\le k\le I_n+A_n$, and $n$ large
enough. Also,\vspace*{1pt}
on the event $\{T^*>I_n+1\}$, we have $\B(0,r_{I_n}-h_{I_n})=A^*_{I_n}(N)$
after the $I_n$th wave. Thus,
\[
W_{I_n}(\Sigma_{I_n})= \bigl|\B(0,n)\bigr|-\bigl|\B
(0,r_{I_n}-h_{I_n})\bigr|
\le2c_1 A_n n^{d-1}\times h_{I_n}.
\]
A key feature of the flashing process is that
explorers stopped, at time $I_n$, outside
$\B(z,3h_{I_n})\cap\Sigma_{I_n}$ cannot settle in $\tilde\C(z)$.
In other words, knowing~$\GG_k$, the covering of a family
of cells $\{\tilde\C(z_j), j=1,\ldots,\NN\}$ are independent
events if $\|z_i-z_j\|\ge6h_{I_n}$ for $i\not= j$.
Now, there is an integer $\NN$ and sites $\{z_j, j=1,\ldots,\NN\}$
with
\[
\forall i\not= j \qquad\|z_i-z_j\|\ge6h_{I_n}\quad
\mbox{and}\quad\sum_{j\le\NN} \bigl|\B(z_j,6
h_{I_n})\cap\Sigma_{I_n}\bigr|\ge|\Sigma_{I_n}|.
\]
We then get using (\ref{def-constants}),
%
%
\begin{equation}
\label{optimal2} \NN h_{I_n}^d\ge\tfrac{1}{2}c_2
h_{I_n} n^{d-1}.
\end{equation}
Let
\[
\Gamma=\bigl\{j\in[1,\NN]\dvtx W_{I_n}\bigl(B(z_j,
3h_{I_n}) \cap\Sigma_{I_n}\bigr)\le c_3
a_n h_{I_n}^d\bigr\} \quad\mbox{and}\quad
\Gamma^c=[1,\NN]\bs\Gamma.
\]
On $\{T^*>I_n+1\}$,
\begin{eqnarray*}
2c_1 A_n n^{d-1}\times h_{I_n}
&\ge&
W_{I_n}(\Sigma_{I_n})\ge\sum_{j\in\Gamma^c}
W_{I_n}\bigl(B(z_j, 3h_{I_n})\cap
\Sigma_{I_n}\bigr)\\
&\ge&\bigl|\Gamma^c\bigr|\times\bigl(
{c_3 a_n h_{I_n}^d} \bigr).
\end{eqnarray*}
Thus, using the definition of $A_n$ in (\ref{def-aA}), and bound
(\ref{optimal2}) on $\NN$, we obtain
\[
\bigl|\Gamma^c\bigr|\le\frac{2c_1 A_n h_{I_n} n^{d-1}}{c_3 a_n
h_{I_n}^d}\le\frac
{c_1 c_2 c_3 a_n}{2c_3c_1 a_n}
\frac{\NN}{c_2}=\frac{\NN}{2}.
\]
In other words, we have that $|\Gamma|\ge\NN/2$. Now, as already
noticed,
knowing $\GG_{I_n}$, for any subset $I\subset[1,\NN]$,
the events $\{
\tilde{\cal C}(z_j) \subset{\cal A}^*_{I_n+1}(N), j\in I\}$
are independent. By conditioning on $\GG_{I_n}$, we obtain
for $h_0$ large enough,
%
%
\begin{eqnarray}
\label{opti-cond}
&&
P\bigl(\bigl\{T^*>I_n+1\bigr\}\bigr) \nonumber\\
&&\qquad= E \biggl[ {
\sum_{I\subset[1,\NN], |I|\ge\NN/2}\ind_{\Gamma
=I}\times P \bigl( {
\forall i\in I, \tilde{\cal C}(z_j) \subset{\cal A}^*_{I_n+1}(N)
|\GG_{I_n}} \bigr)} \biggr]
\nonumber\\[-8pt]\\[-8pt]
&&\qquad= E \biggl[ {\sum_{I\subset[1,\NN], |I|\ge\NN/2}\ind_{\Gamma
=I}\times
\prod_{i\in I} P \bigl( {\tilde{\cal
C}(z_j) \subset{\cal A}^*_{I_n+1}(N) |\GG_{I_n}}
\bigr)} \biggr]
\nonumber
\\
&&\qquad\le\sup_{z_j\in\Sigma_{I_n}} P \bigl( {{\cal
A}^*_{I_n+1}(N)\supset
\tilde{
\cal C}(z_j), W_{I_n}\bigl(B(z_j,
3h_{I_n})\cap\Sigma_{I_n}\bigr)\le c_3
h_{I_n}^da_n} \bigr)^{\NN/2}.\nonumber\hspace*{-25pt}
\end{eqnarray}
Considering the probability appearing on the right-hand side of
(\ref{opti-cond}), we can think of a coupon-collector problem, where
an album of size $|\tilde{\cal C}(z)|$ has to be filled
when we collect no more than $c_3 h_{I_n}^da_n$ coupons.
Using inequality (\ref{coupon-main}) of
Lem\-ma~\ref{lem-coupon} below, we show that
\[
P \bigl( { \bigl\{ {T^*>I_n+1} \bigr\}} \bigr)\le\exp\biggl(-
\frac{\alpha_1}{4}a_n^2\frac
{c_2}{2}h_{I_n}^{1-d}n^{d-1}
\biggr).
\]
This concludes the proof.

The result about filling an album, that we just mentioned,
is based on the following simple
coupon-collector lemma (together with Proposition~\ref{prop-uniform}),
which we did not find in the vast literature on such problems.
%
%
\begin{itlemma}\label{lem-coupon}
Consider an album of $L$ items
for which are bought independent
random cou\-pons, each of them covering one (or possibly none)
of the possible $L$ items.
If $Y_i$ is the item
associated with the $i$th coupons,
we assume that for positive constants $\alpha_1,\alpha_2$, such that
for any $j=1,\ldots,L$,
%
%
\begin{equation}
\label{coupon1} \frac{\alpha_1}{L} \le P(Y_i=j)\le
\frac{\alpha_2}{L}.
\end{equation}
Let $\tau_L$ be the number of coupons needed to complete the album.
Then, for any $0<A<\frac{1}{4\alpha_2} \log(L)$, we have
%
%
\begin{equation}
\label{coupon-main} P(\tau_L<A L)\le\exp\biggl(-\frac{\alpha_1^2
A^2e^{-2\alpha_2 A}}{4}
\sqrt{ L} \biggr) \le\exp\biggl( {-\frac{\alpha_1^2 A^2}{4}} \biggr).
\end{equation}
\end{itlemma}
\begin{pf}
We denote by $\sigma_i$ the time needed to collect the $i$th
distinct item after having collected $i-1$
distinct items.
The sequence $\{\sigma_1,\sigma_2,\ldots,\sigma_L\}$
is not independent, but if $\YY_k=\sigma(\{Y_1,\ldots,Y_k\})$,
and $\tau(k)=\sigma_1+\cdots+\sigma_k$, then for $i=1,\ldots,L$,
%
%
\begin{equation}
\label{coupon2}\hspace*{28pt} \biggl(1-\frac{\alpha_1(L-i+1)}{L} \biggr
)^k\ge P(
\sigma_{i}>k | \YY_{\tau(i-1)})\ge\biggl(1-\frac{\alpha_2(L-i+1)}{L}
\biggr)^k.
\end{equation}
Indeed, calling ${\cal E}(i-1)$ the set of the first $i-1$ collected items,
%
%
\begin{eqnarray}
\label{coupon3}\quad
P (\sigma_{i}>k | \YY_{\tau(i-1)} )&=&P \bigl(
\{Y_{\tau(i-1)+1},\ldots, Y_{\tau(i-1)+k}\}\subset{\cal E}(i-1) |
\YY_{\tau(i-1)} \bigr)
\nonumber
\\
&=& \bigl(P \bigl( Y \in{\cal E}(i-1) | \YY_{\tau(i-1)} \bigr)
\bigr)^k
\\
&=& \bigl( 1-P \bigl( Y \notin{\cal E}(i-1) | \YY_{\tau
(i-1)} \bigr)
\bigr)^k.\nonumber
\end{eqnarray}
Using (\ref{coupon1}) we deduce (\ref{coupon2}) from (\ref{coupon3}).
Formula (\ref{coupon2}) gives that
\[
\frac{L}{\alpha_1(L-i+1)}\ge E[\sigma_{i}| \YY_{\tau(i-1)}]\ge
\frac{L}{\alpha_2(L-i+1)}
\]
as well as
%
%
\begin{equation}
\label{coupon11} E \bigl[\sigma_{i}^2 |
\YY_{\tau(i-1)} \bigr]\le2\frac
{L^2}{\alpha_1^2(L-i+1)^2}.
\end{equation}
Now, we look for $B\leq\sqrt L$ such that
%
%
\begin{equation}
\label{coupon5} \sum_{i={\sqrt L}}^{B{\sqrt L}} E[
\sigma_{L-i}]\ge2 AL.
\end{equation}
Note that
\[
\sum_{i={\sqrt L}}^{B{\sqrt L}} E[\sigma_{L-i}]
\ge\frac{L}{\alpha_2}\sum_{i={\sqrt L}}^{B{\sqrt L}}
\frac{1}{i+1}\ge\frac{L}{\alpha_2} \log(B).
\]
Thus, condition (\ref{coupon5}) holds for $B\ge\exp(2\alpha_2 A)$,
but recall that $B\le\sqrt L$ also, and this gives a bound on $A$.
Finally, note that
\[
\max\bigl\{ E [\sigma_{L-i} | \YY_{\tau(L-i-1)} ], i={\sqrt L},
\ldots,B{\sqrt L} \bigr\}\le\frac{\sqrt L}{\alpha_1}
\]
and set
\[
X_i=\frac{E [\sigma_{L-i} | \YY_{\tau(L-i-1)} ]
-\sigma_{L-i}}{ ( \sqrt L/\alpha_1 )}\le1.
\]
For $x\le1$, note that $e^x\le1+x+x^2$ to obtain for $0\le\lambda
\le
1$, by successive conditioning,
%
%
\begin{eqnarray}
\label{coupon7} P \Biggl(\sum_{i={\sqrt L}}^{B{\sqrt L}}
\sigma_{L-i}\le AL \Biggr) &\le& P \Biggl( \sum
_{i={\sqrt L}}^{B{\sqrt L}} X_i\ge\alpha_1
A \sqrt L \Biggr)
\nonumber
\\
&\le& e^{-\lambda\alpha_1 A\sqrt{ L}}\prod_{i={\sqrt L}}^{B{\sqrt L}}
\bigl(1+\lambda^2 \sup E \bigl[X_i^2 |
\YY_{\tau(L-i-1)} \bigr] \bigr)
\\
&\le& \exp\biggl(-\lambda\alpha_1 A\sqrt L+\lambda^2
\sum_i \sup E \bigl[X_i^2
| \YY_{\tau(L-i-1)} \bigr] \biggr).
\nonumber
\end{eqnarray}
Finally,
we have, using (\ref{coupon11}),
\[
\sum_{i={\sqrt L}}^{B{\sqrt L}} \sup E
\bigl[X_i^2 | \YY_{\tau
(L-i-1)} \bigr]\le\sum
_{i={\sqrt L}}^{B{\sqrt L}} \alpha_1^2\sup
\frac{E [\sigma^2_{L-i} |\YY_{\tau(L-i-1)} ]}{L} \le2 B \sqrt L.
\]
The results follows as we optimize on $\lambda\le1$ in the upper bound
in (\ref{coupon7}).
\end{pf}

\section{Potential theory estimates}\label{sec-harmonic}
We collect in this section three technical results.
In Corollary~\ref{cor-gauss}, we estimate the difference between the
expected number of independent random walks exiting a ball
$\B(0,n)$ at a distinguished site, whether the random walks
are initially on the origin or are spread over a sphere $\B(0,r_n)$
with $r_n<n$. Corollary~\ref{cor-gauss} is used to bound the
mean number of explorers exiting some large ball from a given
site, and its proof relies on a discrete mean value property
Theorem~\ref{theo-gauss}, which in turns relies on
Blach\`ere's Proposition~\ref{prop-L} written in the \hyperref[app]{Appendix}.
Then, Lemma~\ref{lem-alex} improves an estimate of Lawler, Bramson and
Griffeath in~\cite{lawler92}, dealing with the exit site
distribution from a sphere when the initial position is not
the origin. Indeed, Lemma 5(b) of~\cite{lawler92}, states
that when $d\ge2$, there is a positive constant $J_d$ such that for
any $r>0$, $z\in\B(0,r)$ and $z^*\in\partial\B(0,r)$, we have
%
%
\begin{equation}
\label{hitting-lawler} \P_z\bigl(S(H_r) = z^*\bigr)
\leq\frac{J_d}{(\|z^*\|-\|z\|)^{d-1}}.
\end{equation}
Thus, when $\|z^*\|-\|z\|$ is small, (\ref{hitting-lawler}) is useless.
Since we need bounds on the sum of squares of $\P_z(S(H_r) = z^*)$
over $z\in\B(0,r-h)$ of order $\log(r)$ in $d=2$, and of order
$1/h^{d-2}$ when $d>2$, we establish the following.
%
%
\begin{itlemma}\label{lem-alex}
There is a positive constant $\kappa_G$ such that,
for all $r>0$, if $z\in{\mathbb B}_r$, $z^*\in\partial{\mathbb B}_r$,
then
%
%
\begin{equation}
\label{hitting-alex} \P_z\bigl(S(H_r) = z^*\bigr) \leq
\frac{\kappa_G}{\|z-z^*\|^{d-1}}.
\end{equation}
\end{itlemma}
Finally, we prove the \textit{uniform hitting
property},
Proposition~\ref{prop-uniform}, for
the boundary of a \textit{cell}. Though the property is natural,
the nonspherical nature of a cell, makes its proof tedious.

\subsection{A discrete mean value theorem}\label{sec-gauss}
The following result has interest on its own.
%
%
\begin{theorem}\label{theo-gauss} There are positive constants $K_0$
and $K_a$ such that
for any sequence $\{\Delta_n,n\in\N\}$ with $K_0\le\Delta_n\le n^{1/3}$,
for any $z\in\B_n$ with $n-\|z\|\le1$, we have, setting
$r_n=n-\Delta_n$,
%
%
\begin{equation}
\label{gauss-estimate} \biggl||\B_{r_n}|\times G_n(0,z) - \sum
_{y\in\B_{r_n}} G_n(y,z) \biggr|\le K_a.
\end{equation}
\end{theorem}
%
%
\begin{itremark}
\label{rem-gauss} Note that a related (but
distinct) property was also at
the heart of~\cite{lawler92}. Namely, for $\varepsilon>0$, and
$n$ large enough, if $z\in\B_n$, and $n-\|z\|\ge\varepsilon n$,
%
%
\begin{equation}
\label{heart-lawler} |\B_{n}|\times G_n(0,z) \ge\sum
_{y\in\B_n} G_n(y,z).
\end{equation}
\end{itremark}
We start with proving the following useful corollary of
Theorem~\ref{theo-gauss}.
%
%
\begin{itcorollary}\label{cor-gauss}
In the setting of Theorem~\ref{theo-gauss}, and for
any $\Lambda\subset\partial\B_n$,
%
%
\begin{equation}
\label{prereq-main} \bigl|E \bigl[ {M\bigl(|\B_{r_n}|\id_0,n,
\Lambda\bigr)} \bigr] -E \bigl[ {M(\B_{r_n},n,\Lambda)} \bigr] \bigr|
\le
K_a |\Lambda|.
\end{equation}
\end{itcorollary}
\begin{pf}
Note that (\ref{prereq-main}) holds if for any $z^*\in\partial\B_n$,
%
%
\begin{equation}
\label{prereq1} \biggl||\B_{r_n}|\times\P_0 \bigl(
{S(H_n)=z^*} \bigr) -\sum_{y\in\B_{r_n}}
\P_y \bigl( {S(H_n)=z^*} \bigr) \biggr| \le K_a.
\end{equation}
By a classical decomposition
(Lemma 6.3.6 of~\cite{lawler-limic}), we have
for a finite subset $B\subset\Z^d$, $y\in B$, and $z^*\in\partial B$
%
%
\begin{equation}
\label{har-5} \P_{y} \bigl(S \bigl(H(\partial B) \bigr)= z^* \bigr) =
\frac{1}{2d}\sum_{z\in B,z\sim z^*} G_B(y,z).
\end{equation}
For $B=\B(0,n)$, we replace in (\ref{prereq1}) the value
of $\P_{y}(S(H(\partial B))=z^*)$ by the right-hand side in (\ref{har-5}),
and are left with proving that for any $z\in\B_n$ with $n-\|z\|\le1$,
we have (\ref{gauss-estimate}).
\end{pf}
\begin{pf*}{Proof of Theorem~\ref{theo-gauss}}
When $d\ge3$, we express $G_n(0,z)$ in
term of Green's function (Proposition 4.6.2(a) of~\cite{lawler-limic}),
\[
G_n(0,z)=G(0,z)-\E_z\bigl[G\bigl(0,S(H_n)
\bigr)\bigr].
\]
Now, using Green's function asymptotics (\ref{main-green}),
there is a constant $K_1$ (independent on $n$) such that
%
%
\begin{equation}
\label{prereq5}\qquad\biggl\llvert v_d G_n(0,z)-2
\frac{\alpha(z)}{n^{d-1}}\biggr\rrvert\le\frac{K_1}{n^d}\qquad
\mbox{where }
\alpha(z)=
\E_z \bigl[ {\bigl\|S(H_n)\bigr\|-\|z\|} \bigr].
\end{equation}
In $d=2$, $G_n$ is expressed in terms of the potential kernel
(Proposition~4.6.2(b) of~\cite{lawler-limic})
\[
G_n(0,z)=-a(0,z)+\E_z\bigl[a\bigl(0,S(H_n)
\bigr)\bigr].
\]
Using (\ref{main-potential}), we have
\[
\pi G_n(0,z)=2\alpha(z)/n+O\bigl(1/n^2\bigr).
\]

Now, $r_n^d=n^d-d\Delta_n n^{d-1}+O(\Delta_n^2 n^{d-2})$, so that
using (\ref{prereq5}), and the hypothesis $\Delta_n=O(n^{1/3})$, and
$0\le n-\|z\|\le1$
%
%
\begin{eqnarray}
\label{prereq6} |\B_{r_n} | G_n(0,z)
&=& \bigl(
{r_n^d+O\bigl(r_n^{d-1}\bigr)}
\bigr) \biggl( {2\frac{\alpha(z)}{n^{d-1}}+O\biggl(\frac{1}{n^d}\biggr
)} \biggr)
\nonumber
\\
&=& \bigl( {n^d-d\Delta_n n^{d-1}+O\bigl(
\Delta_n^2 n^{d-2}\bigr)+O\bigl(n^{d-1}
\bigr)} \bigr)\nonumber\\[-8pt]\\[-8pt]
&&{}\times\biggl( {2\frac{\alpha(z)}{n^{d-1}}+O\biggl(\frac
{1}{n^d}\biggr)}
\biggr)
\nonumber
\\
&=& 2\alpha(z) (n-d\Delta_n)+O(1).\nonumber
\end{eqnarray}
Since $\{\|S_n\|^2-n, n\in\N\}$ is a martingale (with the natural
filtration),
\[
\E_z\bigl[\bigl\|S(H_n)\bigr\|^2\bigr]-\|z
\|^2] =\E_z[H_n]= \sum
_{y\in\B_n} G_n(y,z).
\]
Using $n-\|z\|\le1$, this yields for a constant $K_l$,
%
%
\begin{equation}
\label{prereq7} \biggl\llvert\sum_{y\in\B_n}
G_n(y,z)-2\alpha(z)n\biggr\rrvert\le K_l.
\end{equation}
We now invoke Proposition~\ref{prop-L} of the \hyperref[app]{Appendix}.
There is $K_b$ such that for $z\in\B_n$ with $n-\|z\|\le1$,
%
%
\begin{eqnarray}
\label{prereq8} \biggl\llvert\sum_{y\in\A(r_n,n)}
G_n(y,z)-2\alpha_0(z)d\Delta_n \biggr
\rrvert\le K_b\nonumber\\[-8pt]\\[-8pt]
&&\eqntext{\mbox{where } \alpha_0(z)=
\E_z \bigl[\bigl\|S(H_n)\bigr\|-\|z\| | H_n<H(B_{r_n})
\bigr].}
\end{eqnarray}
From (\ref{prereq7}) and (\ref{prereq8}), we obtain
%
%
\begin{equation}
\label{prereq9} \biggl\llvert\sum_{y\in\B_{r_n}}
G_n(y,z) -2n\alpha(z)+2\alpha_0(z)d
\Delta_n\biggr\rrvert\le K_l+K_b.
\end{equation}
Now, from (\ref{prereq6}) and (\ref{prereq9}) we obtain
for a constant $K_2$,
\[
\biggl\llvert\biggl( {|\B_{r_n}| G_n(0,z)-\sum
_{y\in\B_{r_n}} G_n(y,z)} \biggr)- 2 \bigl( {
\alpha_0(z)-\alpha(z) } \bigr) d\Delta_n\biggr\rrvert
\le K_2.
\]
Now, from
\begin{eqnarray*}
\bigl|\alpha_0(z)-\alpha(z)\bigr| &\le&\P_z \bigl(
{H(B_{r_n})<H_n} \bigr)\\
&&{}\times\bigl( {
\alpha_0(z) +\E_z \bigl[ {\bigl\|S(H_n)\bigr\|-\|z\|
| H_n>H(B_{r_n})} \bigr]} \bigr),
\end{eqnarray*}
and the Gambler's ruin estimate,
for $K_0>0$ and $z\in\AA(n-1,n)$,
\[
\P_z \bigl( {H(B_{r_n})<H_n} \bigr)\le
\frac{K_0}{\Delta_n},
\]
we deduce that
\[
\Delta_n\bigl|\alpha_0(z)-\alpha(z)\bigr|\le2
\Delta_n \P_z \bigl( {H(B_{r_n})<H_n}
\bigr)\le2K_0.
\]
The desired result follows.
\end{pf*}

\subsection{\texorpdfstring{Proof of Lemma \protect\ref{lem-alex}}{Proof of Lemma 5.1}}
We follow\vspace*{1pt} the proof of Lemma 5(b) of~\cite{lawler92}.
Set $D:=\|z-z^*\|$.
Let $O'$ be a closest point to
$(1+ \frac{D}{4r})z^*$ in $\B(z^*, \frac{D}{4})$.
We define $B_1':= \B(O',\frac{D}{4})$, $B_2':= \B(O',\frac{D}{2})$,
and we note that
\[
\bigl\|z-z^*\bigr\| \leq\bigl\|z-O'\bigr\|
\]
and, for all $x$ in $\partial B_2'$, the triangle inequality
$\|z-z^*\| \leq\|z-x\| + \|x- z^*\|$ implies that
\[
{\min_{x\in\partial B_2'}}\|z-x\| \geq{D \over3}.
\]
Now, define
\[
\tau:= \inf\bigl\{ t> 0 \dvtx S(t) \in\{z\}\cup B_r^c
\bigr\} \quad\mbox{and}\quad\tau':= \inf\bigl\{ t> 0 \dvtx S(t) \in
B_1'\cup\partial B_2' \bigr\}.
\]
By a last exit decomposition, together with the strong Markov property,
%
%
\begin{eqnarray}
\label{alex-11} \P_z\bigl(S(H_r) = z^*\bigr)
&=&
G_r(z,z)\P_{z^*}\bigl(S(\tau) = z\bigr)
\nonumber
\\
&\leq& G_r(z,z) \P_{z^*}\bigl(S\bigl(\tau'
\bigr)\in B_2'\bigr) \max_{x\in\partial B'_2}
\P_x\bigl(S(\tau) = z\bigr)
\nonumber\\[-8pt]\\[-8pt]
&=& \P_{z^*}\bigl(S\bigl(\tau'\bigr)\in
B_2'\bigr) \max_{x\in\partial B_2'} G_r(x,z)
\nonumber
\\
&\leq& \P_{z^*}\bigl(S\bigl(\tau'\bigr)\in
B_2'\bigr) \max_{x\in\partial B_2'} G_{r+D}(x,
z).\nonumber
\end{eqnarray}
A Gambler's ruin estimate yields, for some positive constant $c$,
\[
\P_{z^*}\bigl(S\bigl(\tau'\bigr)\in B_2'
\bigr) \leq{c\over D}.
\]
The desired result follows from (\ref{alex-11}) and the previous bound,
after we show that for a constant $c$,
such that for all $x$ satisfying $\|x-z\| \geq{D\over3}$,
%
%
\begin{equation}
G_{r+D}(x,z) \leq{c\over D^{d-2}}.
\end{equation}
On the set $V:= \B(z,{D\over4})$, the map $y\mapsto G(x,y)$ is harmonic.
By Harnack's inequality, we have
%
%
\begin{equation}
\label{alex-13} G_{r+D}(x,z) \leq{c\over D^d}\sum
_{y\in V} G_{r+D}(x,y) = {c\over D^d}
\E_x[Y],
\end{equation}
where $c$ is a positive constant, and
$Y$ is the number of visits of $V$ before time~$H_{r+D}$.
By taking the supremum over the entering site of $V$ in (\ref{alex-13}),
\[
G_{r+D}(x,z) \leq{c\over D^d}\sup_{y\in V}
\E_y[Y].
\]
It remains to show that $\sup_{y\in V}\E_y[Y]\le J D^2$, for some
positive constant $J$. This is identical to (2.10) of~\cite{lawler92},
and we omit this last step.

\subsection{\texorpdfstring{Proof of Proposition \protect\ref{prop-uniform}}{Proof of Proposition 3.1}}
\label{sec-uniform}

For $j\geq0$, consider $z_j$ in $\Sigma_j$.
We show that for positive constants $\alpha_1$, $\alpha_2$,
and for all $z^*$ in ${\cal C}(z_j)$, we have (\ref{uniform-inequality}).
The random walk has initial condition $S(0)=z_j$.

First, when $z^*=z_j$, $S(\sigma_j)=z_j$ if and only if $X_j=1$. This
happens with probability $1/h_j^d$, and gives the result in this case.

Assume $z^*\in\C(z_j)\bs\{z_j\}$. We recall that
the unbiased Bernoulli variable $Y_j$ decides whether the explorer
can flash upon exiting either a sphere
or an annulus. More precisely, we draw
$R_j$ with density $g_j$ given in (\ref{density}), and if $Y_j=1$
(resp., $Y_j=0$) the walk flashes upon exiting the ball of center $z_j$
and radius $R_j\wedge(r_j+h_j-\|z_j\|)$ [resp., $\AA(r_j-R_j,r_j+R_j)$]
provided $S(\sigma_j)\in\C(z_j)$.

\textit{Step} 1: \textit{Flashing when exiting a sphere $(Y_j=1)$.}
We first prove the upper bound when $Y_j=1$ and $X_j=0$.
It is obvious that
\[
z^*\in\partial\B\bigl(z_j,\bigl\|z^*-z_j\bigr\|\bigr)\qquad
\mbox{but } z^*\notin\partial\B\bigl(z_j,\bigl\|z^*-z_j\bigr\|-1
\bigr).
\]
Thus, $R_j\in\ ]\|z^*-z_j\|-1,\|z^*-z_j\|]$,
and there is a constant $C$ such that
%
%
\begin{eqnarray}
\label{har-6}
&&
P \bigl( {X_j=0, Y_j=1,R_j
\in\, \bigl]\bigl\|z^*-z_j\bigr\|-1,\bigl\|z^*-z_j\bigr\|\bigr]} \bigr
)\nonumber\\[-8pt]\\[-8pt]
&&\qquad\le C
\frac{\|z^*-z_j\|^{d-1}}{h_j^d}.\nonumber
\end{eqnarray}
On the other hand, by
(\ref{hitting-position}) of Section~\ref{sec-green},
%
%
\begin{eqnarray}
\label{har-4}
&&
P \bigl( {S(\sigma_j)=z^* | X_j=0,
Y_j=1,R_j\in\, \bigl]\bigl\|z^*-z_j\bigr\| -1,\bigl\|
z^*-z_j\bigr\|\bigr]} \bigr)\nonumber\\[-9pt]\\[-9pt]
&&\qquad\le\frac{c_2}{\|z^*-z_j\|^{d-1}}.\nonumber
\end{eqnarray}
The upper bound in the case $\{X_j=0, Y_j=1\}$ follows from
(\ref{har-6}) and (\ref{har-4}).

We now turn to the lower bound when $Y_j=1$ and $X_j=0$.
Since we want a lower bound,
we consider the event that the walk flashes on $z^*$
when exiting a sphere
only in the case where $ |\|z^*\| - r_j | < h_j /2$.
Note that
by Lemma~\ref{lem-border}, $z^*$ has a nearest neighbor, say $z$,
which satisfies
\[
\|z-z_j\|\le\bigl\|z^*-z_j\bigr\|-\frac{1}{4\sqrt{d}}.
\]
This means that if
$h\in V:=[\|z^*-z_j\|-1/(4\sqrt{d}),\|z^*-z_j\|[$, then
$z^*\in\partial\B(z_j, h)$. Thus
\begin{eqnarray*}
P_{z_j} \bigl( {S(\sigma_j)=z^*} \bigr)
&\ge&
P(X_j=0, Y_j=1, R_j\in V)\times
\inf_{h\in V} \P_{z_j} \bigl( {S\bigl(H\bigl(\partial
\B(z_j,h)\bigr)\bigr)=z^*} \bigr)
\\[-2pt]
&\ge& \frac{ch^{d-1}}{h_j^d} \inf_{h\in V} \P_{z_j} \bigl( {S\bigl(H
\bigl(\partial\B(z_j,h)\bigr)\bigr)=z^*} \bigr)
\\[-2pt]
&\ge& \frac{ch^{d-1}}{h_j^d}\times\frac{c_1}{h^{d-1}} \qquad\mbox{[using
(\ref{hitting-position})]}.
\end{eqnarray*}
The lower bound in the case $Y_j=1, X_j=0$, and
$ |\|z^*\| - r_j | < h_j /2$ is obtained.

\textit{Step} 2: \textit{Flashing when exiting an annulus $(Y_j=0)$.}
The upper bound for this case is close to the case $Y_j=1$.
It is obvious that
%
%
\begin{eqnarray}
z^*\in\partial\AA\bigl(r_j - \bigl|\bigl\|z^*\bigr\|-r_j\bigr|,
r_j +\bigl|\bigl\|z^*\bigr\|-r_j\bigr|\bigr)\nonumber\\[-2pt]
&&\eqntext{\mbox{but } z^*\notin\partial
\AA\bigl(r_j - \bigl|\bigl\|z^*\bigr\|-r_j\bigr|+1, r_j +\bigl|\bigl\|z^*
\bigr\|-r_j\bigr|-1\bigr).}
\end{eqnarray}
Thus necessarily, $R_j\in\ ]|\|z^*\|-r_j|-1,|\|z^*\|-r_j|]$, and
\[
P \bigl( {Y_j=0,R_j\in\, \bigl]\bigl|\bigl\|z^*\bigr\|-r_j\bigr|-1,\bigl|
\bigl\|z^*\bigr\|-r_j\bigr|\bigr]} \bigr) \le C\frac{|\|z^*\|
-r_j|^{d-1}}{h_j^d}.
\]
For $h>0$, define $\D_h=\AA(r_j-h,r_j+h)$.
It is enough to prove that for some constant $c$, and for any $h$
such that $z^*\in\partial\D_h$ (and $h\in\ ]|\|z^*\|-r_j|-1,\break|\|z^*\|
-r_j|]$),
%
%
\begin{equation}
\label{right-alex} \P_{z_j} \bigl( {S\bigl(H\bigl(\D_h^c
\bigr)\bigr)=z^*} \bigr)\le\frac{c}{h^{d-1}}.
\end{equation}
Note the following fact.
If $\|z^*\|>\|z_j\|$, and the walk exits $\D_h$ at $z^*$, then
the walk exits $\B(0,r_j+h)$ at $z^*$,
whereas if $\|z^*\|<\|z_j\|$, and the walk exits $\D_h$ at $z^*$, then
the walk enters $\B(0,r_j-h)$ at $z^*$. In both cases,
Lemma 5(b) of~\cite{lawler92} yields (\ref{right-alex}).
(Actually, Lemma 5(b) of~\cite{lawler92} is formulated to cover
only the case $\|z^*\|>\|z_j\|$, but its proof covers both cases.)\vadjust{\goodbreak}

We turn now to the lower bound.
By Lemma~\ref{lem-border}, $z^*$ has a nearest neighbor, say~$z$,
%
%
\begin{equation}
\label{simple-nn} \bigl|\bigl\|z\bigr\|-r_j \bigr|\le\bigl|\bigl\|
z^*\bigr\|-r_j \bigr|-
\frac{1}{4\sqrt{d}}.
\end{equation}
This means that if
$h\in V:=[ |\|z^*\|-r_j |-1/(4\sqrt{d}), |\|z^*\|-r_j|[$, then
\mbox{$z^*\in\partial\D_h$}.

We only need to consider the case $ |\|z^*\| - r_j | \geq h_j /2$.
It is enough to prove, for $h\in V$, $z^*\in\C(z_j)\cap\partial\D_h$,
and for some constant $c$ (that depends on $d$), that
%
%
\begin{equation}
\label{goal3} \P_{z_j} \bigl(S\bigl(H\bigl(\D_h^c
\bigr)\bigr)= z^* \bigr) \geq\frac{c}{h^{d-1}}.
\end{equation}
Let $y^*$ be the closest site of $\partial\B(0,r_j)$ to the segment
$[0,z^*]$,
and let $x^*$ be in $\R^d$ given by
\[
x^*=\biggl(r_j+\frac{h}{2}\biggr)\frac{z^*}{\|z^*\|}.
\]
Note that if $\tilde z\sim z$ and $\tilde z$ satisfies
(\ref{simple-nn}), then $\|x^*-\tilde z\|<\inf_{z\in\partial\D_h}
\|x^*-z\|$,
and we define
\[
R^*=\frac{1}{2} \Bigl( {\bigl\|x^*-\tilde z\bigr\|+\inf_{z\in\partial\D
_h} \bigl\|
x^*-z\bigr\|}
\Bigr).
\]
Define
\[
\tilde\D_h=\AA\biggl(r_j-\frac{h}{2},r_j+
\frac{h}{2}\biggr)\quad\mbox{and set}\quad\Gamma= \B\bigl(x^*,
R^*\bigr)\cap
\partial\bigl(\tilde\D_h^c\bigr).
\]
Thus, if $\|z^*\|>\|z_j\|$, then $\Gamma$ is the boundary of
the lower hemisphere of the ball $\B(x^*,R^*)$. We need also
the time
$\tau^+=\inf\{n\ge1\dvtx S(n)\in\D_h^c\cup\{z_j\}\}$.
By a last exit decomposition, and the strong Markov property,
we have
%
%
\begin{eqnarray}
\label{har-1} \qquad\P_{z_j} \bigl(S\bigl(H(\partial\D_h)\bigr)
= z^* \bigr)
&=& G_{\D_h} (z_j, z_j)
\P_{z^*} \bigl(S\bigl(\tau^+\bigr) = z_j \bigr)
\nonumber
\\
&\ge& G_{\D_h} (z_j, z_j) \P_{z^*}
\bigl(H(\Gamma) < \tau^+ \bigr) \min_{x\in\Gamma}\P_{x} \bigl(S(
\tau) = z_j \bigr)
\\
&\ge& \P_{z^*} \bigl(H(\Gamma) < \tau^+ \bigr) \min_{x\in\Gamma}G_{\D_h}
(x, z_j).\nonumber
\end{eqnarray}
Since $z^*\in{\cal C}(z_j)$, we have $\|y^*-z_j\|\le h_j/2$, so that
$y^*$ and $z_j$
can be connected by 10 overlapping balls of radius $h_j / 10$
in such a way that, applying Harnack's inequality 10 times
(see Theorem 6.3.9 in~\cite{lawler-limic})
to the harmonic map $y\mapsto G_{\D_h}(x,y)$,
we can estimate from below the last factor
in (\ref{har-1}). For any $x\in\Gamma$,
\[
G_{\D_h} (x, z_j) \ge c_H^{10}
G_{\D_h} \bigl(x, y^*\bigr).
\]
We use again Harnack's inequality on
the harmonic functions $x\mapsto G_{\D_h}(x,y^*)$, to obtain
\[
\min_{x\in\Gamma}G_{\D_h} \bigl(x, y^*\bigr) \ge c_H
G_{\D_h} \bigl(x', y^*\bigr),\vadjust{\goodbreak}
\]
where $x'\in\B(x^*,R^*/2)$ and $\|x'-y^*\|\in[\frac{h}{4}-1,\frac{h}{4}]$.
The purpose of choosing $x'$ is to have $y^*\in\B(x',h/4)$, and
$\B(x',h/2)\subset\D_h$ so that
$G_{\D_h}(x',y^*)\ge G_{\B(x',h/2)}(x',y^*)$.

When dimension is 2, the classical
expansion of $G_{\B(x',h/2)}(x',\cdot)$ (see Proposition 6.3.5 of
\cite{lawler-limic}) gives with a constant $K_2$,
%
%
\begin{equation}
\label{green-d2}
G_{\B(x',h/2)}\bigl(x',y^*\bigr)\ge
\frac{2}{\pi} \log\biggl( {\frac{h/2}{\|
x'-y^*\|}} \biggr)-\frac{K_2}{\|x'-y^*\|}
\ge\frac{2}{\pi} \log2- \frac{4K_2}{h}.\hspace*{-28pt}
\end{equation}
When $h$ is large enough, $G_{\B(x',h/2)}(x',y^*)\ge\log(2)/\pi$.

When dimension is larger than $2$, by using
(\ref{main-green}), there is a constant $K_d$ such that,
when $h$ is large enough,
%
%
\begin{eqnarray}
\label{har-2} G_{\D_h} \bigl(x', y^*\bigr)
&=& G
\bigl(x', y^*\bigr)-\E_{x'} \bigl[G\bigl(S\bigl(H\bigl(
\D_h^c\bigr)\bigr), y^*\bigr) \bigr]
\nonumber\\[-8pt]\\[-8pt]
&\ge& \frac{C_d}{h^{d-2}} \bigl( {4^{d-2}-1} \bigr)- \frac{K_g}{h_j^{d}}
\ge\frac{K_d}{h^{d-2}}.\nonumber
\end{eqnarray}
As a consequence of (\ref{har-2}), we just need to prove
that the first factor in (\ref{har-1}) is of order $1/h$, at least.
We realize the event $\{ H(\Gamma)<\tau^+\}$ in two moves: the walk first
hits the sphere $\B(x^*,R^*/2)$, and then exits from the cap
$\partial\B(x^*,R^*)\cap\tilde\D_h$,
%
%
\begin{eqnarray}
\label{flashing-18} \P_{z^*} \bigl( {H(\Gamma)<\tau^+} \bigr)
&\ge&
\frac{1}{2d}\P_{\tilde z} \bigl( {H\bigl(\B\bigl(x^*,R^*/2\bigr)
\bigr)<H\bigl(\B^c\bigl(x^*,R^*\bigr)\bigr)} \bigr)
\nonumber\\[-8pt]\\[-8pt]
&&{}\times\inf_{y\in\partial\B(x^*,R^*/2)} \P_y \bigl( {S \bigl(H \bigl
( {
\B^c\bigl(x^*,R^*\bigr)} \bigr) \bigr)\in\tilde\D_h}
\bigr).\nonumber
\end{eqnarray}
The first factor in the right-hand side of (\ref{flashing-18}) is of
order $1/R^*$, that is, of order $1/h$. To deal with the second factor, we
invoke Harnack's inequality to have for $y\in\partial B(x^*,R^*/2)$, and
for $x''$ the closest point of $\Z^d$ to $x^*$,
%
%
\begin{equation}
\label{flashing-19}\quad\P_y \bigl( {S \bigl(H \bigl( {\B^c
\bigl(x^*,R^*\bigr)} \bigr) \bigr)\in\tilde\D_h} \bigr) \ge
c_H \P_{x''} \bigl( {S \bigl(H \bigl( {\B^c
\bigl(x^*,R^*\bigr)} \bigr) \bigr)\in\tilde\D_h} \bigr).
\end{equation}
We invoke now (\ref{hitting-position}) to obtain for some constant $K_3$,
%
%
\begin{equation}
\label{flashing-20} \P_{x''} \bigl( {S \bigl(H \bigl( {\B^c
\bigl(x^*,R^*\bigr)} \bigr) \bigr)\in\tilde\D_h} \bigr) \ge
c_1\frac{| \partial\B(x^*,R^*)\cap\tilde\D_h|}{|\partial\B
(x^*,R^*)|}\ge K_3.
\end{equation}
We gather (\ref{flashing-18}),
(\ref{flashing-19}) and (\ref{flashing-20}) to obtain the desired
lower bound.

%
%
\begin{appendix}\label{app}
\section{\texorpdfstring{Coupling and proof of Lemma \lowercase{\protect\ref{theo-coupling}}}
{Coupling and proof of Lemma 1.3}}
We give a precise definition of our coupling.
To avoid heavy notation, we write the coupling algorithm
as a pseudo-code.

First of all, we draw $N$ independent sequences
of independent Bernoulli and continuous random variables
$((X_{k,l},Y_{k,l},R_{k,l}\dvtx l\geq0)\dvtx1\leq k\leq N)$
as in Section~\ref{sec-flashing}.
In addition, we call $(U_k\dvtx k\geq1)$
the sequence of the increments
of a generic independent simple\vadjust{\goodbreak} random walk on ${\mathbb Z}^d$.
From these two sources of randomness we extract
our explorer and flashing explorer trajectories
with their associated clusters.
The flashing times will be adapted to the flashing explorer trajectories
as in Section~\ref{sec-flashing}.

Our pseudo-code is made of two loops of size $N$ that make precise
the previous description.
With the first loop we build our $N$ random walk trajectories
$((S_i(t)\dvtx0\leq t\leq\tau_i)\dvtx1\leq i\leq N)$ with their
associated clusters $A(1),\ldots,A(N)$.
Step by step,
within this first loop,
we also define pieces of the flashing explorers
trajectories $S_1^*,\ldots,S_N^*$.
With the second loop we complete the trajectories
of the flashing explorers to build the associated cluster $A^*(N)$.
During the algorithm, $t_k \in{\mathbb N}$
stands for the time up to which the trajectory of flashing explorer $k$
has been defined ($k\in\{1;\ldots;N\}$).
We use the same $t$ for the time governing the evolution of each
simple random walk $S_i$.
The index $j$ is updated before adding each random walk increment
to the partial sum of $S_i$ and $S_j^*$. The updating procedure
uses the index $j'$ described in Section~\ref{sec-coupling},
and we denote by $\Delta= U$
the increment. Each encountered $U$ stands for the first unused
random variable in the sequence $(U_k\dvtx k\geq1)$.\looseness=-1

The main advantage of the pseudo-code formalism
is that it allows, through the assignment operator ``$\leftarrow$,''
expressions of the kind $j\leftarrow\max(j,j')$ or $t_j\leftarrow t_j+1$
rather than $j(\theta+ 1) = \max(j(\theta),j'(\theta))$
and $t_{j(\theta+1)}(\theta+1) = t_{j(\theta)}(\theta) + 1$
with $\theta$ a discrete
parameter ordering the sequence of our elementary moves.
It makes also implicit identities like
$t_k(\theta+1) = t_k(\theta)$
for any quantity $t_k$ that does not need to be updated.
Our following pseudo-code can be re-written
in a classical inductive way with $\theta$
running through
$\{(i,t)\in\{1;\ldots;N\}\times{\mathbb N}\dvtx t\leq\tau_i\}$
according to lexicographic order.
Marks $^{(a)}$ and $^{(b)}$ refer
to remarks (a) and (b) below:
\begin{eqnarray*}
&& A(0)\leftarrow\varnothing;
\\
&& \mbox{For $i=1$ to $N$}
\\
&&\cases{ %
j\leftarrow i;
\vspace*{1pt}\cr
t\leftarrow0;\qquad t_j \leftarrow0;
\vspace*{1pt}\cr
S_i(t)\leftarrow0;\qquad S_j^*(t_j) \leftarrow0%
\hfill\mbox{$\bigl[$\textit{Note that
$S_i(t) = S_j^*(t_j)$}$\bigr]$}
\vspace*{2pt}\cr
\mbox{While $S_i(t)\in A(i-1)$}
\vspace*{2pt}\cr
\cases{\mbox{If
$t_j$ is a flashing time %
for explorer
$j^{(a)}$, then}
\vspace*{2pt}\cr
\cases{j'\leftarrow%
\begin{array} {l} \mbox{unique index$^{(b)}$ $k\in\{1;\ldots;i\}
\setminus
\{j\}$}
\\[2pt]
\mbox{such that $S_k^*(t_k) = S_j^*(t_j)
= S_i(t)$};
\end{array}
\cr
\mbox{If $t_{j'}$
is not a flashing time %
for explorer $j'^{(a)}$,
then } j\leftarrow j';
\vspace*{1pt}\cr
\mbox{Otherwise }%
j\leftarrow
\max\bigl(j,j'\bigr); } %
\vspace*{2pt}\cr
\Delta\leftarrow U;
S_i(t+1)\leftarrow S_i(t) + \Delta;
S_j^*(t_j+1)\leftarrow S_j^*(t_j)
+ \Delta
\cr
\hfill\mbox{$\bigl[$\textit{so that $S_i(t+1) =
S_j^*(t_j+1)$}$\bigr]$}
\cr
t\leftarrow t+1;
t_j\leftarrow t_j +1;}
\cr
A(i)\leftarrow A(i-1)\cup
\bigl\{S_i(t)\bigr\};
\cr
i\leftarrow i+1;}
\\
&& A^*(0)\leftarrow\varnothing;
\\
&& \mbox{For $k=1$ to $N$}
\\
&&\cases{ %
\mbox{While $t_k$ is not a flashing time
for explorer $k^{(a)}$ or $S_k^*(t_k)
\in A^*_{k-1}$}
\vspace*{2pt}\cr
\cases{S_k^*(t_k + 1)
\leftarrow S_k^*(t_k) + U;
\vspace*{2pt}\cr
t_k\leftarrow
t_k +1;}
\vspace*{2pt}\cr
A^*(k)\leftarrow A^*(k-1)\cup\bigl\{S_k^*(t_k)
\bigr\};
\vspace*{2pt}\cr
k\leftarrow k+1;}
\end{eqnarray*}

\textit{Remarks}:
\begin{longlist}[(a)]
\item[(a)]
Recall that for $l=j,j'$ or $k$,
$S_l^*$ is defined up to time $t_l$
as well as its associated flashing times.
\item[(b)]
One checks by induction on $i$ that just after the instruction
``$A(i)\leftarrow A(i-1)\cup\{S_i(t)\}$,'' we have
%
%
\begin{equation}
\label{afterwhile} A(i) = \bigl\{S_1^*(t_1);
\ldots;S_i^*(t_i)\bigr\} \quad\mbox{and}\quad \bigl|A(i)\bigr| = i.
\end{equation}
To do so, one checks by induction on $t<\tau_i$,
that
%
%
\begin{eqnarray}
A(i-1) &=& \bigl\{S_1^*(t_1);\ldots;S_{j-1}^*(t_{j-1});
S_{j+1}^*(t_{j+1});\ldots;S_{i}^*(t_i)
\bigr\} \quad\mbox{and}\nonumber\\[-8pt]\\[-8pt]
\quad\bigl|A(i-1)\bigr| &=& i -1.\nonumber
\end{eqnarray}
Since we always have $S_j^*(t_j)= S_i(t)$,
this proves by induction that $j'$ is well defined.
\end{longlist}

The key observation is that for each increment $U$,
the index of the explorer
that follows this increment
depends on the whole previous construction,
but the value of $U$ does not depend on it.
As a consequence, we build independent random walks
$S_1,\ldots,S_N$
coupled with independent flashing random walks
$S_1^*,\ldots,S_N^*$.
Then, one simply checks by induction on $i$ and $k$ that
%
%
\begin{equation}
\label{rightcluster}\quad A(i) = \bigl\{S_1(
\tau_1);\ldots;S_i(\tau_i)\bigr\}
\quad\mbox{and}\quad A^*(k) = \bigl\{S_1^*\bigl(\tau_1^*\bigr);
\ldots;S_k^*\bigl(\tau_k^*\bigr)\bigr\}
\end{equation}
for all $1\leq i,k\leq N$.

Finally, define
$(\bar t_1,\ldots,\bar t_N)$ and $(\tau_1^*,\ldots,\tau_N^*)$
the values of $(t_1,\ldots,t_N)$ at the end of
the first and last cycle, respectively.
Since $t_1,\ldots,t_N$ can only increase during our loops,
we have $\tau_k^*\geq\bar t_k$ for all $k$.
Then (\ref{afterloop1}) and (\ref{afterloop2})
follow from (\ref{afterwhile}) and (\ref{rightcluster}).

\section{Time spent in an annulus (by Blach\`ere)}\label{sec-blachere}

This section is devoted to an asymptotic expansion of the expected
time spent in an annulus $\A(r_n,n)$ for $r_n<n$,
when the random walk is started at some point $z$ within the annulus,
and before it exits the outer shell.
%
%
\begin{itproposition}\label{prop-L} There are positive constants
$K_0,K_b$, such that for any
sequence $\{\Delta_n,n\in\N\}$ with $K_0\le\Delta_n\le n^{1/3}$,\vadjust{\goodbreak}
for any $z \in\A(r_n,n)$, we have setting $r_n=n-\Delta_n$,
%
%
\begin{equation}
\label{seb-main} \biggl\llvert\sum_{y\in\A(r_n,n)}
G_n(z,y) - \bigl( {2d\Delta_n\alpha_0(z)-d\bigl(n-
\|z\|\bigr)^2} \bigr)\biggr\rrvert\le K_b \bigl( {\bigl(n-\|z\| \bigr)
\vee1} \bigr)\hspace*{-35pt}
\end{equation}
with
\[
\alpha_0(z)=\E_z \bigl[ {\bigl\|S(H_n)\bigr\|-\|z\|
|H \bigl( {B^c(0,n)} \bigr)<H \bigl( {B(0,r_n)} \bigr)}
\bigr].
\]
\end{itproposition}
\begin{pf}
Our strategy is to decompose a path into successive strands lying
entirely in the annulus. The first strand is special since the
starting point is any $z\in\A(r_n,n)$. The other strands, if any,
start all on $\partial\B(0,r_n)$. We estimate the time spent inside
the annulus for each strand. Let us remark that we make use of
three facts: (i) precise asymptotics for Green's function, (ii)
$(G(0,S(n)), n\in\N)$ is a martingale and (iii)
$(\|S(n)\|^2-n, n\in\N)$ is a martingale.

Choose $z\in\A(r_n,n)$.
We define the following stopping times $(D_i,U_i, i\ge0)$,
corresponding to the ${i}$th downward and
upward crossings of the sphere of radius~$r_n$. Let
$\theta(n)$ act on trajectories by time-translation of $n$-units.
Let $\tau=H(B_{r_n})\wedge H_n$, $D_0=U_0=0$ and
\[
D_1=\tau\id_{H(B_{r_n})<H_n}+\infty\id_{H_n<H(B_{r_n})}.
\]
If $D_1<\infty$, then $U_1 =H_{r_n}\circ\theta(D_1)+D_1$, whereas
if $D_1=\infty$, then we set $U_1=\infty$.
We now proceed by induction, and assume $D_i,U_i$ are defined.
If $D_i=\infty$, then $D_{i+1}=\infty$, whereas if $D_i<\infty$
(and necessarily $U_i<\infty$), then
\[
D_{i+1}=U_i+ ( {\tau\id_{\tau=H(B_{r_n})}+\infty
\id_{\tau
=H_n}} )\circ\theta(U_i)
\]
and
\[
U_{i+1}
=D_{i+1}+H_{r_n}\circ\theta(D_{i+1}).
\]
With this notation, we can write
%
%
\begin{eqnarray}
\label{tps} \sum_{y\in\A(r_n,n)} G_n(z,y) & = &
\E_z [ {\tau} ]+ \sum_{i=1}^{\infty}
\E_z \bigl[ {\tau\circ\theta(U_{i})\id_{D_i<\infty}}
\bigr]
\nonumber\\[-8pt]\\[-8pt]
&=&\E_z[\tau]+\P_z ( {D_1<\infty} )\times
\I(z),\nonumber
\end{eqnarray}
where
%
%
\begin{equation}
\label{def-I}\quad \I(z)=\sum_{i=1}^\infty
\E_z \bigl[ {\tau\circ\theta(U_i) | D_i<
\infty} \bigr] \prod_{j=1}^{i-1} \bigl( {1-
\P_z (D_{j+1}=\infty| D_{j}<\infty)} \bigr).
\end{equation}
Now, we compute each term of the right-hand side of (\ref{tps}).

We have divided the proof into three steps.

\textit{Step} 1:
First, we show that there is a positive constant $K$ (independent of
$z$ and~$n$) such that when $
z\in\A(r_n,n)$, then
%
%
\begin{equation}
\label{D1bulk} \biggl\llvert\P_z ( {D_1<\infty} )-
\frac{\alpha_0(z)}{\Delta_n} \biggr\rrvert\le\frac{K}{\Delta_n^2}
\bigl( {\bigl(n-\|z\|\bigr)\vee1}
\bigr).
\end{equation}
Note that when $z\in\B(0,n)$, and $n-\|z\|\le1$, (\ref{D1bulk}) yields
%
%
\begin{equation}
\label{D1outer} \biggl\llvert\P_z ( {D_1<\infty} )-
\frac{\E_z [ {\|S(\tau)\|-\|z\| | D_1=\infty}
]}{\Delta_n} \biggr\rrvert\le\frac{K}{\Delta_n^2}.
\end{equation}
Second, we show that for $z\in\A(r_n,n)$, and $i\ge1$,
%
%
\begin{eqnarray}
\label{D1inner}
&&
\biggl\llvert\P_z ( {D_{i+1}=\infty|
D_i<\infty} )\nonumber\\
&&\quad{}- \frac{\E_z [ { ( {\|S(U_i)\|-\|S(D_{i+1})\|} ) \id
_{D_1\circ\theta(U_i)<\infty} | D_{i}<\infty} ]}{\Delta_n}\biggr\rrvert
\\
&&\qquad\le
\frac{K}{\Delta_n^2}.\nonumber
\end{eqnarray}
Our starting point is the classical Gambler ruin estimate, which
in dimension~2, reads with the potential kernel instead of
Green's function,
%
%
\begin{equation}
\label{gambler-eq}\qquad \P_z ( {D_1<\infty} ) =
\frac{G(0,z)-\E_z [ {G(0,S(\tau)) | D_1=\infty}
]}{
\E_z [ {G(0,S(\tau)) | D_1<\infty} ]
-\E_z [ {G(0,S(\tau)) | D_1=\infty} ]}.
\end{equation}
We now expand Green's function (resp., the potential kernel)
using asymptotics (\ref{main-green}) [resp., (\ref{main-potential})].
For this purpose, it is convenient to define a random variable
\[
X(z)=\frac{1}{\|z\|} \bigl( {\bigl\|S(\tau)\bigr\|^2- \|z\|^2}
\bigr).
\]
Note that for any $z\in\A(r_n,n)$, $X(z)/\|z\|$ is small. Indeed,
%
%
\begin{equation}
\label{ann-2} \frac{X(z)}{\|z\|}=\frac{ ( {\|S(\tau)\|-\|z\|} ) (
{\|S(\tau)\| +\|z\|} )} {
\|z\|^2}.
\end{equation}
Since $\Delta_n=n-r_n=O(n^{1/3})$, we have for $n$ large enough,
%
%
\begin{equation}
\label{ann-3} \frac{X(z)}{\|z\|}\le\frac{2(n+1)\Delta_n}{(n-\Delta
_n)^2} \le\frac{8 \Delta_n}{n}
\quad\!\mbox{and}\quad\!\sup_{z\in\A(r_n,n)} \biggl( {\frac{|X(z)|}{\|z\|}}
\biggr)^3\le\frac{8^3 \Delta_n^3}{n}\times\frac{1}{n^2}.\hspace*{-35pt}
\end{equation}
More precisely, $X(z)$ is of order
$2(\|S(\tau)\|-\|z\|)$. Indeed,
$\Delta_n^3\le n$, and (\ref{ann-2}) yields
%
%
\begin{eqnarray}
\label{Xz-general}
&& X(z)=2 \bigl( {\bigl\|S(\tau)\bigr\|-\|z\|} \bigr)+ \biggl( {
\frac{(\|S(\tau)\|-\|z\|)^2}{\|z\|}} \biggr)\nonumber\\[-8pt]\\[-8pt]
&&\quad\Longrightarrow\quad\bigl| X(z)-2
\bigl( {\bigl\|S(\tau)\bigr\|-\|z\|}
\bigr) \bigr|\le\frac
{1}{\Delta_n}.\nonumber
\end{eqnarray}
When dimension $d>2$, we set $\eta(d)=\frac{d-2}{2}$. In
order to use Green's function asymptotics (\ref{main-green}),
we express $S(\tau)$ in terms of $X(z)$ as follows:
%
%
\begin{equation}
\label{ann-1} \frac{1}{\|S(\tau)\|^{d-2}}=\frac{1}{\|z\|^{d-2}} \biggl( {1+
\frac{X(z)}{\|z\|}} \biggr)^{-\eta(d)}.
\end{equation}
We have a constant $K_d$ such that
%
%
\begin{eqnarray}
\label{ann-4}
&&
\biggl\llvert\biggl( {1+\frac{X(z)}{\|z\|}} \biggr
)^{-\eta(d)}\nonumber\\
&&\quad{}-
\biggl( {1-\eta(d)\frac{X(z)}{\|z\|}+ \eta(d)\frac{\eta(d)+1}{2} \biggl
( {
\frac{X(z)}{\|z\|}} \biggr)^2} \biggr)\biggr\rrvert\\
&&\qquad\le
\frac{K_d}{n^2}.\nonumber
\end{eqnarray}
For $d>2$ and any $z\not= 0$,
(\ref{main-green}), (\ref{ann-1}) and (\ref{ann-4}) yield
%
%
\begin{equation}
\label{ann-5} \biggl\llvert G\bigl(0,S(\tau)\bigr)-G(0,z)-\eta(d) C_d
\biggl( {-\frac{X(z)}{\|z\|^{d-1}}+\frac{\eta(d)+1}{2} \frac{X(z)^2}{\|
z\|^{d}}} \biggr)
\biggr\rrvert\le\frac{K_d}{n^d}.\hspace*{-35pt}
\end{equation}
In dimension 2, the potential kernel asymptotic yields for $K_2>0$,
%
%
\begin{equation}
\label{ann-d2} \biggl\llvert a\bigl(0,S(\tau)\bigr)-a(0,z)-\frac{1}{\pi
}
\biggl( {\frac{X(z)}{\|z\|
}+\frac{1}{2} \frac{X(z)^2}{\|z\|^{2}}} \biggr)\biggr
\rrvert\le\frac{K_2}{n^2}.
\end{equation}
In view of (\ref{ann-d2}), we assume
henceforth that (\ref{ann-5}) holds, but in $d=2$,
we think of $\eta(d) C_d=1/\pi$, and $\frac{\eta(d)+1}{2}=1/2$.

Using (\ref{gambler-eq}) and (\ref{ann-5}), we obtain
%
%
\begin{eqnarray}
\label{ann-6}
&&\P_z ( {D_1<\infty} )\nonumber\\
&&\qquad=
\biggl({ \E_z \bigl[ {X(z)
|D_1=\infty} \bigr]-\bar C(z)+
O\biggl(\frac{1}{n}\biggr)}\biggr)\nonumber\\[-8pt]\\[-8pt]
&&\qquad\quad{}\times \biggl(
\E_z \bigl[ {X(z) |D_1=\infty} \bigr]- \E_z \bigl[ {X(z)
|D_1<\infty} \bigr]\nonumber\\
&&\qquad\quad\hspace*{86.5pt}{}+\underline{C}(z)-
\bar C(z)+O\biggl(\frac{1}{n}\biggr)\biggr)^{-1},\nonumber
\end{eqnarray}
where
%
%
\begin{eqnarray}
\label{ann-7} \bar C(z)&=&\frac{\eta(d)+1}{2}\E_z \biggl[ {
\frac{X^2(z)}{\|z\|} \Big| D_1=\infty} \biggr]\quad
\mbox{and}\nonumber\\[-8pt]\\[-8pt]
\underline{C}(z)&=&\frac{\eta(d)+1}{2} \E_z \biggl[ {
\frac{X^2(z)}{\|
z\|} \Big| D_1<\infty} \biggr].\nonumber
\end{eqnarray}
Using (\ref{ann-3}),
we have some rough estimates on $\bar C$ and $\underline{C}$. For any
$z\in\A(r_n,n)$,
%
%
\begin{equation}
\label{C-rough} \bar C(z)=O \biggl( {\frac{\Delta_n^2}{n}} \biggr)= O
\biggl( {
\frac{1}{\Delta_n}} \biggr)\quad \mbox{and}\quad \underline{C}(z)=O \biggl( {
\frac{\Delta_n^2}{n}} \biggr)=O \biggl( {\frac{1}{\Delta_n}} \biggr).\hspace*{-35pt}
\end{equation}
Using (\ref{Xz-general}),
we have better estimates for $\bar C$ and $\underline{C}$.
%
%
\begin{eqnarray}
\label{C-better}
\bar C(z)&=&d\frac{(n-\|z\|)^2}{\|z\|}+O \biggl( {\frac
{(n-\|z\|)\vee
1}{n}}
\biggr),\nonumber\\[-8pt]\\[-8pt]
\underline{C}(z)&=&d\frac{(\|z\|-r_n)^2}{\|z\|}+O \biggl( {\frac
{(\|z\|
- r_n)\vee1}{n}}
\biggr).\nonumber
\end{eqnarray}
The rough estimates (\ref{C-rough}) together with (\ref{Xz-general})
allow us to derive from (\ref{ann-6})
an estimate for $\P_z(D_1<\infty)$, for any $z\in\A(r_n,n)$.
%
%
\begin{eqnarray}
\label{D1-general}
&&\P_z ( {D_1<\infty} )\nonumber\\
&&\qquad=
\frac{ \E_z [ {\|S(\tau)\|
-\|z\| |D_1=\infty} ]
+O ( {{1}/{\Delta_n}} )} {
\E_z [ {\|S(\tau)\|-\|z\| |D_1=\infty} ]-
\E_z [ {\|S(\tau)\|-\|z\| |D_1<\infty} ]+O (
{{1}/{\Delta_n}} )}\hspace*{-12pt}
\\
&&\qquad=\frac{\alpha_0(z)+O ( {{1}/{\Delta_n}} )} {
\Delta_n(1+O ( {{1}/{\Delta_n}} ))}.\nonumber
\end{eqnarray}
This yields (\ref{D1bulk}) since $\alpha_0(z)\le
1+(n-\|z\|)\vee1\le2 (n-\|z\|)\vee1$.

\textit{Case where $z\in\partial B(0,r_n)$.}
On $\{D_1=\infty\}$, we have
%
%
\begin{equation}
\label{ann-11} X(z)=2 \bigl( {\bigl\|S(\tau)\bigr\|-\|z\|} \bigr)+O\biggl(
\frac{1}{\Delta_n}\biggr).
\end{equation}
On $\{D_1<\infty\}$, we have
\[
X(z)=2 \bigl( {\bigl\|S(\tau)\bigr\|-\|z\|} \bigr)+O\biggl(\frac{1}{n}\biggr).
\]
This implies
%
%
\begin{equation}
\label{C-expansion} \bar C(z)=d\frac{\Delta_n^2}{\|z\|}+O\biggl(\frac
{\Delta_n}{n}
\biggr)\quad \mbox{and}\quad \underline{C}(z)=O\biggl(\frac{1}{n}\biggr).
\end{equation}
Thus
%
%
\begin{eqnarray}
\label{ann-13}
\quad&&\P_z ( {D_1=\infty} )\nonumber\\
&&\qquad=
\frac{2\E_z [ {\|z\|-\|
S(\tau)\| | D_1<\infty} ]+\underline{C}(z)+O({1}/{n})} {
\E_z [ {X(z) |D_1=\infty} ]-\E_z [ {X(z)
|D_1<\infty} ]
+\underline{C}(z)-\bar C(z)+
O({1}/{n})}
\nonumber\\[-8pt]\\[-8pt]
&&\qquad= \frac{\E_z [ {\|z\|-\|S(\tau)\| | D_1<\infty}
]+O({1}/{n})} {
\Delta_n+O(1)}
\nonumber
\\
&&\qquad= \frac{\E_z [ {\|z\|-\|S(\tau)\| | D_1<\infty}
]}{\Delta_n}+O\biggl(\frac{1}{\Delta_n^2}\biggr).\nonumber
\end{eqnarray}
In order to obtain (\ref{D1inner}),
we write (\ref{ann-13}) on $\{D_i<\infty\}$, and $z=S(U_i)$
as follows. There is a constant $K$ such that on the event
$\{D_{i}<\infty\}$,
%
%
\begin{equation}
\label{seb-bis1} \biggl\llvert\E_{S(U_i)} [ {\id_{ D_{i+1}=\infty}} ]-
\frac{\E_{S(U_i)} [ { ( {\|S(U_i)\|-\|S(\tau)\|} )
\id_{D_1\circ\theta(U_i)<\infty}} ]}{\Delta_n\times
\P_{S(U_i)} ( {D_1<\infty} )} \biggr\rrvert\le\frac
{K}{\Delta_n^2}.\hspace*{-30pt}
\end{equation}
Note that (\ref{ann-13}) implies that $\P_{S(U_i)} ( {D_1<\infty
} )=1+
O(1/\Delta_n)$, so that (\ref{seb-bis1}) reads as we integrate over
$\{D_{i}<\infty\}$ with respect to $\E_z$
%
%
\begin{eqnarray}
\label{seb-bis2}
&&
\biggl\llvert\P_{z} ( {D_{i+1}=\infty,
D_{i}<\infty} )- \frac{\E_z [ {\id_{D_i<\infty} ( {\|S(U_i)\|-\|S(\tau
)\|
} ) \id_{D_1\circ\theta(U_i)<\infty}} ]}{\Delta_n}\biggr\rrvert\hspace*{-32pt}\nonumber\\[-8pt]\\[-8pt]
&&\qquad\le
\frac{K\P_z(D_i<\infty)}{\Delta_n^2}.\nonumber
\end{eqnarray}
We obtain (\ref{D1inner}) as we divide both sides of (\ref{seb-bis2})
by $\P_z(D_i<\infty)$.

\textit{Step} 2:
We show now that for any $z\in\A(r_n,n)$, we have
%
%
\begin{equation}
\label{tau-bulk} \bigl| \E_z [ {\tau} ]- \bigl( {d \Delta_n
\alpha_0(z)-d\bigl(n-\|z\| \bigr)^2} \bigr) \bigr| \le K \bigl( {\bigl(n-\|z
\|\bigr)\vee1} \bigr).
\end{equation}
When $z\in B_n$ and $n-\|z\|\le1$, (\ref{tau-bulk}) reads
%
%
\begin{equation}
\label{tau-outer} \bigl| \E_z [ {\tau} ]- \bigl( {d \Delta_n
\alpha_0(z)-d\bigl(n-\|z\| \bigr)^2} \bigr)\bigr|\le K.
\end{equation}
When $z\in\A(r_n,n)$, and $i\ge1$, we show that
%
%
\begin{eqnarray}
\label{tau-inner}
&&\biggl\llvert\frac{\E_z [ {\tau\circ\theta(U_i)
|D_i<\infty}
]}{d\Delta_n^2}\nonumber\\
&&\quad\hspace*{0pt}{}- \frac{\E_z [ { ( {\|S(U_i)\|-\|S(D_{i+1})\|} )\id
_{D_1\circ\theta(U_i) <\infty} |D_i<\infty} ]}{\Delta_n} \biggr
\rrvert\\
&&\qquad\le\frac{K}{\Delta_n^2}.\nonumber
\end{eqnarray}
Using that $\{\|S(n)\|^2-n, n\in\N\}$ is a martingale
and the optional sampling theorem (see Lemma 3 of~\cite{lawler92}),
\begin{eqnarray*}
\E_z [ {\tau} ]&=&\E_z \bigl[ {\bigl\|S(\tau)
\bigr\|^2} \bigr]-\| z\|^2=\|z\|\times\E_z \bigl[
{X(z)} \bigr]
\\
&=&\|z\|\times\bigl(\E_z \bigl[ {X(z) |D_1=\infty}
\bigr]\P_z(D_1=\infty)\\
&&\hspace*{32pt}{} + \E_z \bigl[ {X(z)
|D_1<\infty} \bigr]\P_z(D_1<\infty)
\bigr).
\end{eqnarray*}
Thus, using (\ref{ann-6}), simple algebra yields
%
%
\begin{equation}
\label{ann-16} \E_z [ {\tau} ]=\|z\|\times\bigl( {\bigl(\underline
{C}(z)-\bar C(z)\bigr) \P_z(D_1<\infty)+\bar C(z)}
\bigr)+O(1).
\end{equation}
By recalling (\ref{C-better}) and (\ref{D1bulk}),
\begin{eqnarray*}
\E_z [ {\tau} ]&=&d \biggl(\! { \bigl( {\bigl(\|z\|-r_n\bigr)^2-\bigl(n-
\| z\|\bigr)^2+O(\Delta_n)} \bigr) \biggl( {
\frac{\alpha_0(z)}{\Delta_n}+O\biggl(\frac{(n-\|z\|)\vee1}{\Delta
_n^2}\!\biggr)}\! \biggr)} \!\biggr)
\\
&&{}+d\bigl(n-\|z\|\bigr)^2 + O\bigl(\bigl(n-\|z\|\bigr)\vee1\bigr)
\\
&=&d\bigl(2\|z\|-n-r_n\bigr)\alpha_0(z)+d\bigl(n-\|z\|\bigr)^2
+ O\bigl(\bigl(n-\|z\|\bigr)\vee1\bigr)
\\
&=&d\Delta_n \alpha_0(z)-d\bigl(n-\|z\|\bigr)^2+O
\bigl( {\bigl(n-\|z\|\bigr)\vee1} \bigr).
\end{eqnarray*}
This yields (\ref{tau-bulk}).

Assume now that $z\in\partial B(0,r_n)$. From (\ref{ann-16}), we have
\[
\E_z [ {\tau} ]=\|z\|\times\bigl( {\bigl(\bar C(z)-
\underline{C}(z)\bigr) \P_z(D_1=\infty)+
\underline{C}(z)} \bigr)+O(1).
\]
We use (\ref{D1inner}) and (\ref{C-expansion}) to obtain
%
%
\begin{eqnarray}
\label{ann-20}\qquad \E_z [ {\tau} ]
&=& \|z\| \biggl( { \bigl( {d
\Delta_n^2+O(\Delta_n)} \bigr) \biggl( {
\frac{\E_z [ {\|z\|-\|
S(\tau)\| | D_1<\infty} ]}{\Delta_n} +O\biggl(\frac{1}{\Delta_n^2}\biggr
)} \biggr)} \biggr)\nonumber\\
&&{}+O(1)
\\
&=& d\Delta_n \E_z \bigl[ {\|z\|-\bigl\|S(\tau)\bigr\| |
D_1<\infty} \bigr]+O(1).\nonumber
\end{eqnarray}
Now, write (\ref{ann-20}) as follows. There is a constant $K$ such that
for any $z\in\partial\B(0,r_n)$,
%
%
\begin{equation}
\label{ann-21} \biggl\llvert\frac{\E_z [ {\tau} ]}{d\Delta_n^2}- \frac
{\E_z [ {\|z\|-\|S(\tau)\|\id_{D_1<\infty}}
]}{\Delta_n \P_z(D_1<\infty)} \biggr
\rrvert\le\frac{K}{\Delta_n^2}.
\end{equation}
Note that by (\ref{ann-13})
$\Delta_n \P_z(D_1<\infty)=\Delta_n+O(1)$ and $|\|z\|-\|S(\tau)\|
\id_{D_1<\infty}|\le1$, thus
%
%
\begin{equation}
\label{ann-main} \biggl\llvert\frac{\E_z [ {\tau} ]}{d\Delta_n^2}-
\frac{\E_z [ {\|z\|-\|S(\tau)\|\id_{D_1<\infty}}
]}{\Delta_n} \biggr
\rrvert\le\frac{K}{\Delta_n^2}.
\end{equation}
We replace $z$ by $S(U_i)$ in (\ref{ann-main}) under the event $\{
D_i<\infty
\}$ to obtain
%
%
\begin{equation}
\label{ann-22} \biggl\llvert\frac{\E_{S(U_i)} [ {\tau} ]}{d\Delta
_n^2}- \frac{\E_{S(U_i)} [ { ( {\|S(U_i)\|-\|S(D_1\circ\theta
(U_i))\|} ) \id_{D_1\circ\theta(U_i)<\infty}} ]} {
\Delta_n} \biggr
\rrvert\le\frac{K}{\Delta_n^2}.\hspace*{-36pt}
\end{equation}
We multiply both sides of (\ref{ann-22}) by $\id_{D_i<\infty}$,
take the expectation on both side of (\ref{ann-22}) and divide
by $\P_z(D_i<\infty)$ to obtain (\ref{tau-inner}).

\textit{Step} 3: For $i\ge1$, we show the following bounds:
%
%
\begin{eqnarray}
\label{ai} 2\ge\gamma_i\geq\frac{1}{4d\sqrt{d}}\hspace*{110pt}\nonumber\\[-8pt]\\[-8pt]
&&\eqntext{\mbox{where }
\displaystyle \gamma_i=\E_{z} \bigl[ { \bigl( {\bigl\|S(U_i)\bigr\|-
\bigl\|S(D_{i+1})\bigr\|} \bigr)\id_{D_{i+1}<\infty} | D_i<\infty}
\bigr].}
\end{eqnarray}
The upper bound is obvious. For the lower bound, first we restrict to
$\{D_i<\infty\}$, so that $U_i<\infty$.
By Lemma~\ref{lem-border}, $S(U_i)$ has a nearest neighbor $x$,
within $\B(0,r_n)$ such that
$\llVert S(U_i)\rrVert- \llVert x\rrVert\geq1/(2\sqrt{d})$,
and (\ref{ai}) is immediate.

\textit{Step} 4: We show (\ref{seb-main}) using (\ref{tps}).
For $p$ such that $1\leq p \leq\infty$, let
%
%
\begin{equation}
\label{seb-9}\qquad \sigma_p=\sum_{i=1}^{p}
\E_z \bigl[ {\tau\circ\theta(U_i) | D_i<
\infty} \bigr]
\prod_{j=1}^{i-1} \bigl( {1-
\P_z ( {D_{j+1}=\infty|D_j<\infty} )} \bigr).
\end{equation}
Now, (\ref{def-I}) reads
$\I(z)=\lim_{p\to\infty} \sigma_p$
(this is the limit of an increasing sequence).
We establish in this step that, for some constant $\tilde K$, any
integer~$n$,
%
%
\begin{equation}
\label{step4} \lim_{p\to\infty}\biggl\llvert1-\frac{\sigma_p}{d\Delta
_n^2} \biggr
\rrvert\le\frac{\tilde K}{\Delta_n}.
\end{equation}
Once we prove (\ref{step4}), we have all the bounds to
estimate the right-hand side of (\ref{tps}). Indeed, using
(\ref{tau-bulk}), (\ref{D1bulk}) and (\ref{step4}), we have
%
%
\begin{eqnarray}
\label{conclu-appendix}
&&
\E_z[\tau]+\P_z (
{D_1<\infty} )\times\I(z)\nonumber\\
&&\qquad= d\Delta_n
\alpha_0(z)-d\bigl(n-\|z\|\bigr)^2 +O \bigl( {\bigl(n-\|z\|\bigr)\vee1}
\bigr)
\nonumber\\[-8pt]\\[-8pt]
&&\qquad\quad{} + \biggl( {\frac{\alpha_0(z)}{\Delta_n} +O \biggl( {\frac{(n-\|z\|
)\vee1}{\Delta_n^2}} \biggr)}
\biggr) \times\bigl( {d \Delta_n^2+O(
\Delta_n)} \bigr)
\nonumber
\\
&&\qquad= 2d\Delta_n\alpha_0(z)-2d \bigl( {n-\|z\|}
\bigr)^2+O \bigl( {\bigl(n-\| z\|\bigr)\vee1} \bigr).\nonumber
\end{eqnarray}
In order now to prove (\ref{step4}), we
introduce first some shorthand notation. For $p$ and $j$ positive integers,
%
%
\begin{eqnarray}
\label{def-amin}
a_p&=&1-\frac{\sigma_p}{d\Delta_n^2},\qquad \alpha_j=
\P_z (D_{j+1}=\infty|D_{j}<\infty)\quad \mbox{and}\nonumber\\[-8pt]\\[-8pt]
\beta_j&=& \frac{\E_z [ {\tau\circ\theta(U_j) |
D_j<\infty} ]} {
d\Delta_n^2}.\nonumber
\end{eqnarray}
With this notation, (\ref{D1inner}) and (\ref{tau-inner}) read as follows:
%
%
\begin{equation}
\label{step2} \biggl\llvert\alpha_j-\frac{\gamma_j}{\Delta_n}\biggr
\rrvert\le\frac
{K}{\Delta^2_n}\quad\mbox{and}\quad \biggl\llvert\beta_j-
\frac{\gamma_j}{\Delta_n}\biggr\rrvert\le\frac
{K}{\Delta^2_n}\qquad \mbox{so that } |
\alpha_j-\beta_j|\le\frac{2K}{\Delta^2_n}.\hspace*{-35pt}
\end{equation}
Let us rewrite (\ref{seb-9}) as
%
%
\begin{equation}
\label{amin-recursion}\quad a_1 = 1-\beta_1 \quad\mbox{and}\quad
a_p=a_{p-1}-\beta_p \prod
_{j=1}^{p-1} (1-\alpha_j) \qquad\mbox{for } p >
1.
\end{equation}
In order to establish (\ref{step4}), we show by induction that
%
%
\begin{equation}
\label{seb-induction} \Biggl\llvert a_p-\prod
_{j=1}^{p} (1-\alpha_j)\Biggr\rrvert
\le\varepsilon_p
\end{equation}
with for $p>1$,
%
%
\begin{equation}
\label{seb16} \varepsilon_p= \varepsilon_{p-1}+
\frac{2K}{\Delta_n^2} \prod_{j=1}^{p-1} (1-
\alpha_j) \quad\mbox{and}\quad \varepsilon_1=\frac{2K}{\Delta_n^2}.
\end{equation}
Note that it is easy to estimate $\varepsilon_p$ from (\ref{seb16}).
By (\ref{step2}) and for $K_0$ large enough
there is a constant $\kappa_S$ such that
\begin{eqnarray*}
\varepsilon_p &\le& \frac{2K}{\Delta_n^2} \Biggl(1+\sum
_{k=1}^p \exp\Biggl(-\sum
_{j=1}^k \alpha_j \Biggr) \Biggr)
\\
&\le& \frac{2K}{\Delta_n^2} \Biggl( {1+\sum_{k=1}^p
\exp\Biggl( {-\sum_{j=1}^k
\frac{\gamma_j}{2\Delta_n}} \Biggr)} \Biggr) \le\frac
{2K}{\Delta^2_n}\kappa_S
\Delta_n= \frac{2K\kappa_S}{\Delta_n}.
\end{eqnarray*}
Now, by (\ref{step2}), (\ref{seb-induction}) holds for $p=1$,
and we assume it holds for $p-1$.
Then
%
%
\begin{equation}
\label{seb-10}\qquad (1-\beta_p) \prod_{j=1}^{p-1}
(1-\alpha_j)-\varepsilon_{p-1}\le a_p\le(1-
\beta_p) \prod_{j=1}^{p-1} (1-
\alpha_j)+\varepsilon_{p-1}.
\end{equation}
Then by (\ref{step2}), we have (\ref{seb-induction}) with $\varepsilon_p$
satisfying (\ref{seb16}).

Now (\ref{step4}) follows as we notice that step 3 implies,
together with (\ref{step2}) and for $K_0$ large enough, that
\[
\lim_{p\to\infty} \prod_{j=1}^{p} (1-
\alpha_j)=0.
\]
\upqed
\end{pf}
\end{appendix}

\section*{Acknowledgments}

The authors thank the CIRM for a friendly atmosphere
during their stay as part of a \textit{research in pairs} program.


%

\printaddresses

\end{document}